\documentclass[11pt]{amsart}

\usepackage[english]{babel}
\usepackage[T1]{fontenc}
\usepackage[utf8]{inputenc}

\usepackage{mathpazo}
\usepackage{graphicx}
\usepackage{microtype}
\usepackage{enumerate}
\usepackage{fontawesome}

\usepackage{csquotes}
\usepackage[dvipsnames]{xcolor}
\usepackage{url}

\usepackage{breakurl}
\usepackage[breaklinks,bookmarks, colorlinks, citecolor=MidnightBlue!75!black, urlcolor=MidnightBlue!75!black, linkcolor=MidnightBlue!75!black]{hyperref}

\calclayout

\AtBeginDocument{%
   \def\MR#1{}
}

\usepackage{scalerel}

\usepackage{amsmath, amssymb, amsthm, amscd}
\usepackage{commath, mathtools, mathrsfs}
\usepackage{thmtools}
\usepackage{tikz, tikz-cd}
\usetikzlibrary{matrix, arrows, decorations.pathmorphing}

\theoremstyle{plain}

\newtheorem{theorem}{Theorem}
\newtheorem{lemma}[theorem]{Lemma}
\newtheorem{proposition}[theorem]{Proposition}
\newtheorem{corollary}[theorem]{Corollary}

\newtheorem*{theorem*}{Theorem}
\newtheorem*{lemma*}{Lemma}
\newtheorem*{proposition*}{Proposition}
\newtheorem*{corollary*}{Corollary}

\declaretheorem[numbered=no,name=Theorem A]{TA}
\declaretheorem[numbered=no,name=Theorem B]{TB}
\declaretheorem[numbered=no,name=Theorem C]{TC}

\declaretheorem[numbered=no,name=Section Conjecture]{SC}
\declaretheorem[numbered=no,name=Birational Section Conjecture]{BSC}
\declaretheorem[numbered=no,name={Birational Section Conjecture, alternative formulation}]{BSC2}
\declaretheorem[numbered=no,name=Cuspidalization Conjecture]{CC}

\theoremstyle{definition}
\newtheorem{remark}[theorem]{Remark}
\newtheorem{definition}[theorem]{Definition}
\newtheorem{example}[theorem]{Example}

\newtheorem*{conjecture*}{Conjecture}
\newtheorem*{remark*}{Remark}
\newtheorem*{definition*}{Definition}
\newtheorem*{observation*}{Observation}

\newcommand{\on}{\operatorname}

\newcommand{\A}{\mathbb{A}}

\newcommand{\C}{\mathbb{C}}

\newcommand{\gal}{\on{Gal}}

\newcommand{\G}{\mathbb{G}}

\newcommand{\hz}{\hat{\mathbb{Z}}}
\newcommand{\id}{\on{id}}

\newcommand{\mc}{\mathcal}
\newcommand{\mf}{\mathfrak}

\newcommand{\proj}{\on{Proj}}
\newcommand{\Q}{\mathbb{Q}}

\newcommand{\s}{\subseteq}

\newcommand{\spec}{\on{Spec}}

\newcommand{\uaut}{\underline{\on{Aut}}}
\newcommand{\uisom}{\underline{\on{Isom}}}
\newcommand{\upi}{\underline{\pi}}

\newcommand{\xar}{\xrightarrow}
\newcommand{\Z}{\mathbb{Z}}

\renewcommand{\epsilon}{\varepsilon}
\renewcommand{\H}{\on{H}}

\renewcommand{\O}{\mc{O}}
\renewcommand{\P}{\mathbb{P}}
\renewcommand{\phi}{\varphi}
\renewcommand{\projlim}{\varprojlim}

\renewcommand{\tilde}{\widetilde}
\renewcommand{\hat}{\widehat}

\newcommand\cB{\mathcal{B}}

\usepackage{todonotes}

\title{On the birational section conjecture with strong birationality assumptions}
\date{}
\subjclass[2010]{}

\author{Giulio Bresciani}
\address[G. Bresciani]{CRM Ennio de Giorgi, Collegio Puteano, Office 21, Piazza dei Cavalieri 3, 56126 Pisa}
\email{giulio.bresciani@sns.it}
\thanks{The author was partially supported by the DFG Priority Program "Homotopy Theory and Algebraic Geometry" SPP 1786}

\begin{document}

\begin{abstract}
	Let $X$ be a curve over a field $k$ finitely generated over $\mathbb{Q}$ and $t$ an indeterminate. We prove that, if $s$ is a section of $\pi_{1}(X)\to\operatorname{Gal}(k)$ such that the base change $s_{k(t)}$ is birationally liftable, then $s$ comes from geometry. As a consequence we prove that the section conjecture is equivalent to the cuspidalization of all sections over all finitely generated fields.
\end{abstract}

\maketitle
%~ \tableofcontents

\section{Introduction}

Given $X$ a geometrically connected, smooth curve over a field $k$ with absolute Galois group $\gal(k)$, let $\mc{S}_{X/k}$ be the set of sections of $\pi_{1}(X)\to\gal(k)$ modulo the action of $\pi_{1}(X_{\bar{k}})$ by conjugation, these are usually called Galois sections. A Galois section is \emph{geometric} if it is associated with a $k$-rational point of $X$ and \emph{cuspidal} if it is associated with a $k$-rational point in the boundary $\bar{X}\setminus X$, where $\bar{X}$ is the smooth completion of $X$. Grothendieck made the following conjecture, nowadays called the \emph{section conjecture}, in a letter to Faltings \cite{gro97}.

\begin{SC}
    For every smooth, geometrically connected hyperbolic curve $X$ over a field $k$ finitely generated over $\Q$, every Galois section of $X$ is either geometric or cuspidal.
\end{SC}

A birational version of the conjecture, called the \emph{birational section conjecture}, has been studied as it seems more approachable \cite{koe05} \cite{sti15} \cite{saty21}. We can define a set $\mc{S}_{k(X)/k}$ of sections of $\gal(k(X))\to\gal(k)$ analogously to the above, its elements are called birational Galois sections. A birational section is cuspidal if it comes from a rational point of $\bar{X}$.

\begin{BSC}
    For every smooth, geometrically connected curve $X$ over a field $k$ finitely generated over $\Q$, every birational Galois section of $X$ is cuspidal.
\end{BSC}

Recall that a Galois section of $\mc{S}_{X/k}$ is \emph{birationally liftable} if it's in the image of $\mc{S}_{k(X)/k}\to\mc{S}_{X/k}$ \cite[\S 2.1.1]{sti15} (this is slightly different from the definition given in \cite[Definition 266]{sti13}). Here is an alternative formulation of the birational section conjecture, see \autoref{bireq}.

%Clearly, the birational section conjecture implies that every birationally liftable section of $X$ is either geometric or cuspidal. On the other hand, if $s\in \mc{S}_{k(X)/k}$ is a birational section whose image in $\mc{S}_{U/k}$ is geometric or cuspidal for every open subset $U\subset X$, then $s$ is cuspidal since $\mc{S}_{k(X)/k}\simeq\projlim_{U}\mc{S}_{U/k}$ \cite[Lemma 259]{sti13}. Because of this, the birational section conjecture is equivalent to the following.

\begin{BSC2}
    For every smooth, geometrically connected curve $X$ over a field $k$ finitely generated over $\Q$, every birationally liftable Galois section of $X$ is either geometric or cuspidal.
\end{BSC2}

% Hence, the birational section conjecture for all curves over $k$ is equivalent to the statement that every birationally liftable section of every curve over $k$ is geometric or cuspidal.

%The birational section conjecture states that birationally liftable sections are geometric or cuspidal.

In this note we show that, if we strengthen the birationality assumption, we can obtain the desired result. Namely, we show that if $t$ is an indeterminate and $s$ is a Galois section such that the base change $s_{k(t)}$ is birationally liftable, then $s$ is geometric or cuspidal. 

\subsection{Known results}

J. Koenigsmann \cite{koe05} proved that the birational section conjecture holds over finite extensions of $\Q_{p}$. Clearly, one would like to pass from local fields to number fields. Moreover, M. Sa\"idi and M. Tyler \cite{saty21} have proved that the birational section conjecture for number fields implies it for finitely generated extensions of $\Q$.

J. Stix \cite{sti15} obtained partial results about the passage from local fields to number fields, let us describe them. Fix $X$ a smooth, projective curve over a number field $k$ and $s\in\mc{S}_{k(X)/k}$ a Galois section of $k(X)$. Using Koenigsmann's results, $s$ induces a point $x_{\nu}\in X(k_{\nu})$ for every place $\nu$.

Given an open subset $U\s X$, choose some spreading out $\tilde{U}\to\spec \mf{o}_{k,S}$ of $U$, where $\mf{o}_{k}$ is the ring of integers of $k$ and $S$ is some finite set of places. Let $N_{U}$ be the set of places of $k$ such that $x_{\nu}$ is \emph{not} integral: we have that $N_{U}$ depends on the choice of $\tilde{U}$ only up to a finite number of places, hence the Dirichlet density of $N_{U}$ is well defined.

Stix first proves that, if $k$ is a totally real or imaginary quadratic number field, then $N_{U}$ is infinite for some open subset $U\s X$ \cite[Theorem A]{sti15}. Secondly he proves that, if $N_{U}$ has strictly positive Dirichlet density for some open subset $U\s X$, then $s$ comes from a $k$-rational point of $X\setminus U$ \cite[Theorem B]{sti15}. In order to prove the birational section conjecture for totally real or imaginary quadratic number fields, it is then sufficient to bridge the gap between Stix's two results.

\subsection{Our main theorem}

We study the passage from local fields to number fields, too, but we use a different approach. We strengthen the birationality assumption: under this strengthened hypothesis, we obtain a complete result.

\begin{definition*}
    Let $X$ be a curve over a field $k$, $s\in\mc{S}_{X/k}$ a Galois section and $t$ an indeterminate; consider the base change $s_{k(t)}\in\mc{S}_{X_{k(t)}/k(t)}$ of $s$ to $k(t)$ \cite[Definition 27]{sti13}. We say that $s$ is \emph{$t$-birationally liftable} if $s_{k(t)}$ is birationally liftable.
\end{definition*}

\begin{TA}
    \hypertarget{TA}{Let} $X$ be a smooth curve over a field $k$ finitely generated over $\Q$. A Galois section of $X$ is geometric or cuspidal if and only if it is $t$-birationally liftable.
\end{TA}

%~ While the proof of \hyperlink{TA}{Theorem A} is quite technical, the underlying idea is simple. 
A common theme in the study of the section conjecture is the attempt to replicate with Galois sections the geometric constructions we make with points, see e.g. \cite[Chapter 3]{sti13}. The idea underlying the proof of \hyperlink{TA}{Theorem A} is, given $s$ a Galois section of $\P^{1}\setminus\{0,1,\infty\}$, to define ``the polynomial $t-s$'' so that, if $s$ is geometric associated to a $k$-rational point $x \in k\setminus\{0,1\}$, then ``$t-s$'' coincides with the polynomial $t-x$. If $s$ is $t$-birationally liftable then we are able to define this ``polynomial'' by using a birational lifting of $s_{k(t)}$ to produce an element of $\H^{1}(k(t),\hz(1))=\hat{k(t)^{*}}$.

 %~ ``a polynomial with coefficients in non-cuspidal Galois sections of $\P^{1}\setminus\{0,1,\infty\}$'', which conjecturally should correspond to $k\setminus\{0,1\}$. Given $s$ a Galois section of $\P^{1}\setminus\{0,1,\infty\}$, if $s$ is $t$-birational then we are able to define the ``the polynomial $t-s$'' by using a birational lifting of $s_{k(t)}$ to produce an element of $\H^{1}(k(t),\hz(1))=\hat{k(t)^{*}}$ which, if $s$ is geometric associated to a $k$-rational point $x$, coincides with the image of $t-x$ with respect to the natural map $k(t)^{*}\to\hat{k(t)^{*}}$.

This turns out to be sufficient: if we are over a number field $k$ and, for a place $\nu$, the Galois section $s_{k_{\nu}}$ is associated with a point $x_{\nu}\in k_{\nu}\setminus\{0,1\}$, the fact that the ``polynomial $t-s=t-x_{\nu}$'' is defined over $k$ forces $x_{\nu}$ to be $k$-rational and not to depend on $\nu$.

%~ Translating this idea into a proof requires a number of technical tools. One of them is understanding specialization for \emph{ramified} Galois sections.  specialization is done classically only for \emph{unramified} sections, see \cite[Chapter 8]{sti13}. We develop such a theory in \autoref{ramspe}. 
%~ We use a non-standard valuative criterion for proper morphisms of Deligne-Mumford stacks: this allows us to specialize Galois sections just as we specialize points of schemes using the valuative criterion for properness. 

\subsection{Consequences for the section conjecture}

One of the main reasons for studying the birational section conjecture is to reduce the section conjecture to a \emph{lifting}, or \emph{cuspidalization}, problem.

The section conjecture easily implies the following statement, which is sometimes called \emph{cuspidalization conjecture} in the literature.

\begin{CC}
    For every smooth, geometrically connected hyperbolic curve $X$ over a field $k$ finitely generated over $\Q$ and every non-empty open subset $U\subset X$, the map $\mc{S}_{U/k}\to\mc{S}_{X/k}$ is surjective.
\end{CC}

By a limit argument, the cuspidalization conjecture holds if and only if every Galois section of every hyperbolic curve over $k$ is birationally liftable, see \autoref{cuspbl}.

H. Esnault and P.H. Hai showed that the cuspidalization conjecture reduces the section conjecture to the case of $\P^{1}\setminus\{0,1,\infty\}$ \cite[Proposition 7.9]{eh08}. For $\P^{1}\setminus\{0,1,\infty\}$ one may hope to solve the conjecture with explicit computations, such as the $n$-nilpotents obstructions introduced by J. Ellenberg and K. Wickelgren \cite{wic12}. Moreover, the cuspidalization conjecture reduces the section conjecture to the birational section conjecture. Neither the case of $\P^{1}\setminus\{0,1,\infty\}$ nor the birational section conjecture is needed: the section conjecture is equivalent to the cuspidalization conjecture.

\begin{TB}
	\hypertarget{TB}{Let} $X$ be a hyperbolic curve over a field $k$ finitely generated over $\Q$. The following are equivalent.
	\begin{itemize}
		\item For every finitely generated extension $K/k$, every Galois section of $X_{K}$ is either geometric or cuspidal.
		\item For every finitely generated extension $K/k$, every Galois section of $X_{K}$ is birationally liftable.
	\end{itemize}
        As a consequence, the section conjecture is equivalent to the cuspidalization conjecture.
\end{TB}

\subsection{Consequences for the birational section conjecture}

It is well known that it is enough to prove the birational section conjecture for $\P^{1}\setminus\{0,1,\infty\}$. Thanks to \hyperlink{TA}{Theorem A}, it is then enough to show that for every birationally liftable section $s$ of $\P^{1}\setminus\{0,1,\infty\}$ over a number field $k$ the base change $s_{k(t)}$ lifts to $\gal(k(t)(\P^{1}))$. We prove that it is enough to find a much simpler lifting, and we manage to reduce the existence of said lifting to a problem of ``interpolation of Galois sections'', i.e. finding a Galois section with prescribed specializations of some curve over $k(t)$. Let us explain this.

We say that a morphism $X\to V$ is a \emph{family of curves} if it has the form $\bar{X}\setminus D\to V$ where $\bar{X}\to V$ is smooth, proper, with geometrically connected fibers of dimension $1$, and $D\subset\bar{X}$ is a divisor finite étale over $V$. Assume that $V$ is an affine curve over a field of characteristic $0$ and let $F$ be a geometric fiber, then the sequence
\[1 \to \pi_{1}(F)\to\pi_{1}(X)\to\pi_{1}(V)\to 1\]
is exact, see \autoref{families} (we stress that we are assuming $V$ affine, otherwise this is false).

We may define a space $\mc{S}_{X/V}$ of ``global'' sections $\pi_{1}(V)\to\pi_{1}(X)$ modulo conjugation by $\pi_{1}(F)$. Let us call these $\pi_{1}$ sections, in analogy with Galois sections. If $V$ is an affine open subset of $\P^{1}$ over a number field $k$, a $\pi_{1}$ section $r$ of $X\to V$ is uniquely determined by its specializations $r_{v}\in\mc{S}_{X_{v}/k}$, $v\in V(k)$, see \autoref{hilbsp}.
%~ : this follows from the fact that the subgroup generated by the sections $\gal(k)\to\pi_{1}(V)$ associated with rational points of $V$ is dense by Hilbert's irreducibility theorem, see \autoref{hilbgen}.

If $z\in\mc{S}_{k(\P_{1})/k}$ is a birational Galois section of $\P^{1}$, for every open subset $U\subset\P^{1}$ denote by $\iota_{U}(z)$ the image of $z$ in $\mc{S}_{U/k}$. Let $\Delta\subset\P^{1}\times\P^{1}$ be the diagonal.

\begin{TC}
	\hypertarget{TC}{The} following are equivalent.
	\begin{itemize}
		\item The birational section conjecture holds.
		%~ \item For every number field $k$, every open subset $U\s\P^{1}$ and every birational section $s\in\mc{S}_{U/k}$, there exists an open subset $V\s U$ such that the extension $s\times V$ of $s$ to $U\times V\to V$ lifts to $U\times V\setminus\Delta\to V$.
		\item For every number field $k$, every section $z\in\mc{S}_{k(\P_{1})/k}$ and every open subset $U\s\P^{1}$, there exists an open subset $V\s U$ and a $\pi_{1}$ section $r$ of $U\times V\setminus\Delta \to V$ such that the specialization $r_{v}$ is equal to $\iota_{U\setminus\{v\}}(z)$ for every $v\in V(k)$.
		%~ \item For every number field $k$, every open subset $U\s\P^{1}$ and every birational section $s\in\mc{S}_{U/k}$ with lifting $z\in\mc{S}_{k(\P_{1})/k}$, there exists an open subset $V\s U$ and a representative $\gal(k)\to\pi_{1}(U\setminus\{v\})$ of $z|_{U\setminus\{v\}}$ for every $v\in V(k)$ such that the kernel of $\ast_{v}\gal(k)\to\pi_{1}(V)$ is mapped to the identity by the induced homomorphism $\ast_{v}\gal(k)\to\pi_{1}(U\times V\setminus\Delta)$, where $\ast$ denotes the free product in the category of pro-finite groups. 
	\end{itemize}
\end{TC}

Since a $\pi_{1}$ section is uniquely determined by its specializations, the section $r$ in the statement of Theorem C is a lifting to $U\times V \setminus \Delta$ of the base change $(\iota_{U}(z))_{V}\in\mc{S}_{U\times V/V}$.

%~ In short, Theorem C reduces the birational section conjecture to proving that if $z$ is a Galois section of $\gal(k(\P^{1}))\to\gal(k)$ and $U\s\P^{1}$ is an open subset, it is possible to find an open subset $V\s U$ and "compatible" representatives of the sections $z|_{U\setminus v}$, $v\in V(k)$.

\subsection*{Acknowledgements}

I would like to thank an anonymous referee for providing a very large number of comments. This note has become substantially better thanks to his or her work.

\subsection*{Notation and conventions}

Throughout the article, it is tacitly assumed that schemes do not have points of positive characteristic. In particular, fields are of characteristic $0$. Curves are smooth and geometrically connected. If $X$ is a curve, we denote by $\bar{X}$ its smooth completion. The letter $R$ always denotes a DVR with fraction field $K$ and residue field $k$. If $A$ is an abelian group, we denote by $\hat{A}$ the projective limit $\projlim_{n}A/nA$. We write b.l. (resp. $t$-b.l.) as an abbreviation for birationally liftable (resp. $t$-birationally liftable).

\section{Étale fundamental gerbes}

In addition to the classical language of étale fundamental groups, we use the language of \emph{étale fundamental gerbes} \cite[\S 8, \S 9]{bv15}, \cite[Appendix A]{anab}. Étale fundamental gerbes are essentially an alternative point of view on the theory of étale fundamental groups where groups are replaced by a particular type of algebraic stack, i.e. gerbes. Everything that can be done in one language can be translated into the other one and it is actually possible to use both languages at the same time, switching back and forth depending on convenience.

%~ We use both languages, choosing each time which one is convenient depending on the situation. 
Generally speaking, fundamental groups are better for handling geometric arguments about the original varieties, while fundamental gerbes are better for handling arguments about Galois sections. The reason is that the étale fundamental gerbe, in some sense, \emph{is} the space of Galois sections, but constructed directly without passing through fundamental groups.

Furthermore, fundamental gerbes are naturally base point free. This is particularly helpful since we make a lot of specialization arguments involving different fibers of a single family, and doing so while keeping track of base points would be problematic. With fundamental gerbes, the problem just doesn't exist. Moreover, the theory of non-unique specializations in the classical language requires making non-canonical choices, while this is not the case from the point of view of gerbes.

If $X$ is geometrically connected over a field $k$, the étale fundamental gerbe $\Pi_{X/k}$ is a profinite étale stack over $k$ with a morphism $X\to\Pi_{X/k}$ universal among morphisms to finite étale stacks over $k$. The space of Galois sections $\mc{S}_{X/k}$ of $X$ is in natural bijection with the set of isomorphism classes of $\Pi_{X/k}(k)$. Furthermore, the étale fundamental gerbe behaves well with respect to base change in characteristic $0$ \cite[Proposition A.18]{anab}.

Fundamental gerbes are easier to understand when there is a $k$-rational point. Recall that, if $G$ is an affine group scheme over $k$, the classifying stack $\mc{B}_{k}G$ is defined as follows: given a scheme $S$ over $k$, then $\mc{B}_{k}G(S)$ is the groupoid of $G$-torsors for the fpqc topology over $S$, see for instance \cite[Definition 8.1.14]{olsson}. If $G$ is of finite type over $k$ every $G$-torsor is trivialized by an fppf covering $S'\to S$, namely $S'=T$, and hence the fppf topology is sufficient. In general, e.g. if $G$ is pro-finite, we need to use the fpqc topology, see \cite{bv15} for details. If $X$ has a rational point $x\in X(k)$, the Galois group $\gal(\bar{k}/k)$ acts continuously on the étale fundamental group $\pi_{1}(X_{\bar{k}},x)$, this allows us to define an étale fundamental group scheme $\upi_{1}(X,x)$ which is a twisted form of $\pi_{1}(X_{\bar{k}},x)$. The étale fundamental gerbe identifies naturally with the classifying stack $\mathcal{B}_{k}\upi_{1}(X,x)$, and the $k$-rational sections of $\Pi_{X/k}$ correspond to $\upi_{1}(X,x)$-torsors over $k$. 

\subsection{Relative fundamental gerbes}

In order to do specialization arguments, we need to work with families of curves for which there is a short exact sequence of étale fundamental groups. We identify very restrictive assumptions under which this works, these are sufficient for our purposes.

%\begin{lemma}
%    Let $S$ be any connected scheme, and $X\to S$ a family of curves with with geometric fiber $F$. The sequence
%    \[\pi_{1}(F)\to\pi_{1}(X)\to\pi_{1}(S)\to 1\]
%    is exact.
%\end{lemma}

\begin{lemma}\label{families}
	Let $k$ be a field of characteristic $0$ and $A$ a ring which is a localization of a $1$-dimensional regular domain $\tilde{A}$ of finite type over $k$, write $C=\spec A$. If $X\to C$ is a family of curves and $F$ is a geometric fiber, then
	\[1\to\pi_{1}(F)\to\pi_{1}(X)\to\pi_{1}(C)\to 1\]
	is exact.  
\end{lemma}

\begin{proof}
    Assume first that $A$ is a $1$-dimensional regular domain of finite type over $k$.

    If $k=\C$, the sequence of topological fundamental groups is exact since the second homotopy group of $C$ is trivial, and we get an exact sequence of profinite completions by \cite[Proposition 5]{anderson} and \cite[Proposition 3.7]{gjz08}.

    Assume that $k$ is any field of characteristic $0$, it is clearly enough to do the case in which $k$ is algebraically closed. In this case, since everything is of finite type there exists an algebraically closed field $k'$ with embeddings $k'\subset k$ and $k'\subset \C$ such that $X\to C$ descends to a family of curves $X'\to C'$ over $k'$. Since étale fundamental groups are invariant under base change of algebraically closed fields \cite[Exposé X, Corollaire 1.8]{sga1}, the statement follows from the case $k=\C$.

    If $A$ is a localization of a $1$-dimensional regular domain $\tilde{A}$ of finite type over $k$, since $X$ is of finite type over $C=\spec A$ then up to localizing a finite number of elements of $\tilde{A}$ we may assume that $X$ extends to a family of curves $\tilde{X}\to\tilde{C}=\spec\tilde{A}$. We have that $\pi_{1}(C)\simeq\projlim_{U}\pi_{1}(U)$ where $U$ varies among open subsets of $\tilde{C}$ containing $C\subset \tilde{C}$; this follows from the fact that the category of finite étale covers of $C$ is the direct limit of the categories of finite étale covers of $U$ for varying $U$. Analogously, $\pi_{1}(X)\simeq\projlim_{U}\pi_{1}(\tilde{X}_{U})$.

    For every open subset $C\subset U\subset \tilde{C}$, by the preceding case we have a short exact sequence
    \[1\to\pi_{1}(F)\to\pi_{1}(\tilde{X}_{U})\to\pi_{1}(U)\to 1,\]
    hence we obtain a projective system of short exact sequences where the left term is constant. This implies that 
    \[\pi_{1}(X)\simeq\projlim_{U}\pi_{1}(\tilde{X})\times_{\pi_{1}(\tilde{C})}\pi_{1}(U)\simeq \pi_{1}(\tilde{X})\times_{\pi_{1}(\tilde{C})}\pi_{1}(C),\]
    and the statement follows.

\end{proof}

Let $X\to C$ be a family of curves as in \autoref{families}. Define the \emph{relative fundamental gerbe} $\Pi_{X/C}$ by the $2$-cartesian diagram
\[\begin{tikzcd}
	\Pi_{X/C}\rar\dar	& 	\Pi_{X/k}\dar	\\
	C\rar				&	\Pi_{C/k},
\end{tikzcd}\]
there is a structural morphism $X\to\Pi_{X/C}$ over $C$.

%~ If $k'/k$ is a field extension and $c:\spec k'\to C$ is a point, the sequence
%~ \[1\to\pi_{1}(X_{c,\bar{k'}})\to\pi_{1}(X_{\bar{k'}})\to\pi_{1}(C_{\bar{k'}})\to 1\]

If $k'/k$ is a field extension and $c:\spec k'\to C$ is a point, the induced morphism 
\[\Pi_{X_{c}/k'}\to\Pi_{X/C}\times_{C}\spec k'\]
is an isomorphism: using the fact that the étale fundamental gerbe behaves well under base change \cite[Proposition A.18]{anab}, this is a direct consequence of the exactness of the sequence of étale fundamental groups. Because of this, the relative fundamental gerbe is a convenient way of packing the spaces of Galois sections of the fibers without choosing base points. So to speak, it is the ``relative family of spaces of Galois sections''.

\begin{lemma}\label{limit}
	The relative fundamental gerbe $\Pi_{X/C}\to C$ is a projective limit $\projlim \Phi_{i}$ of proper, étale morphisms $\Phi_{i}\to C$ which are gerbes over $C$.
\end{lemma}

\begin{proof}
	Since $\Pi_{X/k}$, $\Pi_{C/k}$ are pro-finite étale, then we may write them as projective limits $\projlim_{i}\Psi_{i}$, $\projlim_{i}\Lambda_{i}$ of finite étale gerbes over $k$ with $\Pi_{X/k}\to\Psi_{i}$, $\Pi_{C/k}\to\Lambda_{i}$ locally full \cite[Definition 3.4, Proposition 3.9]{bv19}. Up to re-indexing, we may assume that $\Pi_{X/k}\to\Pi_{C/k}$ induces a morphism $\Psi_{i}\to\Lambda_{i}$ for every $i$; the fact that the fibers of $X\to C$ are geometrically connected implies that $\Pi_{X/k}\to\Pi_{C/k}$ is locally full (i.e. the homomorphism of geometric fundamental groups is surjective), which in turn implies that $\Psi_{i}\to\Lambda_{i}$ is locally full. The morphism $\Psi_{i}\to\Lambda_{i}$ is a proper étale relative gerbe: it is proper étale since $\Psi_{i},\Lambda_{i}$ are proper étale over $k$, and it is a relative gerbe by \cite[Proposition 3.10]{bv19}. The statement then follows by defining $\Phi_{i}$ as $\Psi_{i}\times_{\Lambda_{i}}C$.
\end{proof}

\section{Non-unique specializations of ramified Galois sections}\label{ramspe}

Classically, the specialization of a Galois section is defined if the original Galois section is \emph{unramified} \cite[\S 8.2]{sti13}, i.e. when the \emph{ramification homomorphism} is trivial. A notion of specialization exists \emph{always}, as long as we don't require that the specialization is unique. If the section is unramified, this "generalized specialization" is unique and coincides with the classical specialization. While a definition of specialization for ramified sections exists in the literature \cite[Corollary 8.9]{sti13}, little else in known beyond the definition itself.

%~ \subsection{Specializing loops}

%~ Let $R$ be a DVR with fraction field $K$ and residue field $k$ and $X\to \spec R$ a proper morphism of schemes. A generic section $\spec K\to X_{K}$ extends uniquely to a section $\spec R\to X$ thanks to the valuative criterion of properness, this gives a specialization $\spec k\to X_{k}$. Suppose now that $X\to\spec R$ is a family of curves, then $\Pi_{X/R}\to\spec R$ is a projective limit of proper morphisms by \autoref{limit}: given a generic Galois section $\spec K\to\Pi_{X_{K}/K}$, we would like to find a specialization. 

%~ The classical valuative criterion for proper morphisms of stacks \cite[Théorème 7.3]{lmb00} does not help us, since it requires passing to an extension of $R$ and thus enlarging the residue field. 

Let $R$ be a DVR with fraction field $K$ and residue field $k$, the \emph{$n$-th root stack} $\sqrt[n]{\spec R}$ of $\spec R$ is defined for every $n$ \cite[Appendix B]{agv08}. If $\pi\in R$ is a uniformizing parameter, then $\sqrt[n]{\spec R}$ is isomorphic to the quotient stack $[\spec R(\sqrt[n]{\pi})/\mu_{n}]$, but the definition of root stack does not depend on the choice of $\pi$. The morphism $\sqrt[n]{\spec R}\to \spec R$ is generically an isomorphism, while the closed fiber with the reduced structure is non-canonically isomorphic to the classifying stack $B\mu_{n}$.

The \emph{infinite root stack} $\sqrt[\infty]{\spec R}$ is the projective limit $\projlim_{n}\sqrt[n]{\spec R}$, see \cite{tv18} for details. Let $c\in\spec R$ be the closed point and $S^{1}_{c}$ the reduced fiber of $\sqrt[\infty]{\spec R}$ over $c$, it is non-canonically isomorphic to the classifying stack $\cB_{k}\hz(1)$ of $\hz(1)$ over $k$. The notation $S^{1}_{c}$ is meant to be reminiscent of the topological space $S^{1}$, which is the classifying space $\mathcal{B}\Z$ of $\Z$. There is a non-canonical isomorphism between the isomorphism classes of $S^{1}_{c}(k)$ and $\H^{1}(k,\hz(1))=\projlim_{n}k^*/k^{*n}=\hat{k^*}$.
 %~ Because of this, the correct topological analogue of $\cB_{k(c)}\hz(1)$ is $\C^{*}$ rather than $S^{1}$, but using the notation $\C^{*}_{c}$ would be extremely confusing.

Recall that a morphism $Y\to X$ of fibered categories over $k$ is \emph{constant} if there exists a factorization $Y\to\spec k\to X$.

\begin{proposition}\label{prorootval}
	Let $X\to C$ be a family of curves as in \autoref{families}, $c\in C$ a closed point with local ring $R$ and $z:\spec K\to\Pi_{X/C}$ a generic Galois section, where $K$ is the fraction field of $C$. Then $z$ extends to a $2$-commutative diagram
	\[\begin{tikzcd}
		\spec K\rar\ar[dr]		&	\sqrt[\infty]{\spec R}\rar[dotted]\dar	&	\Pi_{X/C}\dar	\\
								&	\spec R\rar												&	C.
	\end{tikzcd}\]
	The extension $\sqrt[\infty]{\spec R}\to\Pi_{X/C}$ is unique up to a unique isomorphism. 
	
	We call the induced morphism $S^{1}_{c}\to\Pi_{X_{c}/k(c)}$ the specializing loop $\gamma_{z}(c)$ of $z$ at $c$. An extension $\spec R\to\Pi_{X/C}$ exists if and only if the specializing loop is constant.
\end{proposition}

\begin{proof}
	Using \autoref{limit}, this is a direct consequence of \cite[Theorem 3.1]{giulio-angelo-valuative}.
\end{proof}

\begin{definition}\label{loopdef}
	%~ With notation as in \autoref{prorootval}, the composition $H_{c}\to\sqrt[\infty]{\spec R}\to\Pi_{X/C}$ factorizes as a morphism
	%~ \[H_{c}\to\Pi_{X_{c}/k(c)}\]
	%~ to the fundamental gerbe of the fiber, we call this morphism the \emph{specializing loop} $h_{z}(c)$ of $z$ at $c$. 
	A \emph{specialization} of $z$ at $c$ is any Galois section of $X_{c}/k(c)$ in the essential image of the specializing loop $S^{1}_{c}(k(c))\to\Pi_{X_{c}/k(c)}(k(c))$ of $z$ at $c$ evaluated on $k(c)$. A specialization always exists, but in general it might be not unique. If the specializing loop is constant, then the specialization is unique (the converse is false in general).
\end{definition}

If we fix any section $s\in S^{1}_{c}(k(c))$, the specializing loop $\gamma_{z}:S^{1}_{c}\to\Pi_{X_{c}/k_{c}}$ is constant if and only if the corresponding homomorphism of group schemes $\hz(1)=\uaut(s)\to\uaut(\gamma_{z}(s))$ is trivial. Suppose for simplicity that $k(c)=k$. By the characterization given above, the specializing loop of $z$ at $c$ is constant if and only if the specializing loop of $z_{\bar{k}(C)}$ is constant at $c_{\bar{k}}$: having constant specializing loop in the closed fiber $X_{c}$ only depends on the base change to $\bar{k}$ of the family, nonetheless it encapsulates arithmetic information over $k$, i.e. uniqueness of the specialization.

%~ \begin{lemma}\label{uniqueloop1}
	%~ Let $X\to \spec R$ a geometric fibration, $z\in\Pi_{X_{K}/K}(K)$ a generic section, $k'/k$ an extension of the residue field. Then the specializing loop $\gamma_{z}$ is constant if and only if the base change $\gamma_{z,k'}$ is constant.
	%~ \begin{proof}
		%~ Choose a section of $H_{c}(k)$ so that we may identify $H_{c}=B\hz(1)$, $\Pi_{X_{k}/k}=BG$ for some group scheme $G$ over $k$. The specializing loop is constant if and only if the associated homomorphism $\hz(1)\to G$ is trivial, and this condition is clearly invariant under base change.
	%~ \end{proof}
%~ \end{lemma}

\begin{remark}
	A specialization of a ramified section can be constructed in the language of étale fundamental groups as follows (see also \cite[Corollary 89]{sti13}). Let $R$ be a DVR with fraction and residue fields $K$, $k$ of characteristic $0$, denote by $R^{h}$ a henselianization and $K^{h}$ the fraction field. We may identify $\gal(\bar{k}/k)=\pi_{1}(\spec R^{h})$, and if $X\to \spec R$ is a family of curves then $\pi_{1}(X_{k})\simeq\pi_{1}(X_{R^{h}})$. Since $R^{h}$ is henselian, the Galois group $\gal(\bar{K}/K^{h})$ coincides with the decomposition subgroup and hence it is an extension of $\gal(\bar{k}/k)$ by the inertia subgroup.
	
	Let $\pi\in R^{h}$ be a uniformizing parameter, choose a $n$-th root of $\pi$ for every $n$ in a compatible way and let $K^{h}(\sqrt[\infty]{\pi})$ be the corresponding extension of $K^{h}$, we have that $\gal(\bar{K}/K^{h}(\sqrt[\infty]{\pi}))\simeq\gal(\bar{k}/k)$ defines a section of $\gal(\bar{K}/K^{h})\to\pi_{1}(\spec R^{h})=\gal(\bar{k}/k)$.

	If $\gal(\bar{K}/K)\to\pi_{1}(X_{K})$ is a generic Galois section, it induces a Galois section $\gal(\bar{K}/K^{h})\to\pi_{1}(X_{K^{h}})$ and a specialization is defined by the composition
	\[\gal(\bar{k}/k)\to \gal(\bar{K}/K^{h})\to \pi_{1}(X_{K^{h}})\to\pi_{1}(X_{R^{h}})\simeq\pi_{1}(X_{k}).\]
	Let $I\subset \gal(\bar{K}/K^{h})$ be the inertia subgroup, by composition we have a map $I\to\pi_{1}(X_{\bar{k}})\subset\pi_{1}(X_{k})$ which is called \emph{ramification map}: unramified Galois sections are those for which this map is trivial. 
	
	If we fix a section of $D\to\gal(\bar{k}/k)$, the corresponding specialization induces an action of $\gal(\bar{k}/k)$ on $\pi_{1}(X_{\bar{k}})$, and the ramification map is Galois-equivariant with respect to the natural Galois action on $I\simeq\hz(1)$. The set of specializations is the image in nonabelian Galois cohomology of the ramification map. 
	
	%~ Furthermore, for our arguments it would be necessary to define the specializing loop as a Galois equivariant homomorphism from the inertia group to the étale fundamental group of the fiber: again, this depends on some choices, while the specializing loop as a morphism of gerbes is canonical.
\end{remark}

\begin{example}\label{gmloop}
	Let $R$ be a DVR as in \autoref{families} with fraction field $K$ and valuation $v:K^{*}\to\Z$, this extends naturally to an homomorphism $\hat{K^{*}}\to\hz$ which we still call $v$. The morphism $\G_{m}=\A^{1}_{R}\setminus\{0\}\to \spec R$ is a family of curves.
	
	The section $1\in\G_{m}(R)$ gives an identification $\Pi_{\G_{m}/R}=\cB_{R}\hz(1)$, in particular $\Pi_{\G_{m}/R}(R)=\hat{R^*}$ and $\Pi_{\G_{m,K}/K}(K)=\hat{K^*}$. If $z\in\hat{K^{*}}$ is a generic section, by \autoref{prorootval} the specializing loop $\gamma_{z}$ is constant if and only if $z$ is in the image of $\hat{R^{*}}\to\hat{K^{*}}$, or equivalently if and only if $v(z)=0\in\hz$.
\end{example}

It can be proved that unramified Galois sections in the sense of \cite[Chapter 8]{sti13} are those for which the specializing loop is constant, though we don't need this fact. If the specializing loop is constant, the specialization is clearly unique; this is coherent with the fact that unramified Galois sections have a canonical specialization.

The converse is false in general: if the residue field is algebraically closed, the specialization is unique since there is only one $\hz(1)$-torsor up to equivalence, but the Galois section might be ramified and the specializing loop not constant. For instance, we may choose $R=\C[x]_{(x)}$ in \autoref{gmloop}, the generic Galois section associated with $x\in\hat{\C(x)^{*}}\setminus\hat{\C[x]_{(x)}^{*}}$ has non-constant specializing loop but unique specialization. 
%~ The bottom line is that ramification is a geometric concept, while uniqueness of specialization is an arithmetic one.

Still, in arithmetic situations, it is often the case that uniqueness of specialization forces a constant specializing loop. 

\begin{lemma}\label{uniqueloop}
	Let $T$ be a torus over a field $k$ with a surjective valuation $\nu:k^{*}\to\Z$. Suppose that we have a morphism of gerbes $\gamma:\cB_{k}\hz(1)\to \Pi_{T/k}$. If the image of $\cB_{k}\hz(1)(k)$ in $\Pi_{T/k}(k_{\nu})$ has a finite number of isomorphism classes, then $\gamma$ is constant.
 %, or equivalently the induced homomorphism $f:\hz(1)\to\upi_{1}(T)$ is trivial.
	\begin{proof}
		%~ We may identify $\Pi_{T/k}$ with $\cB_{k} G$ where $G=\projlim_{n}T[n]$ is the Tate module of $T$, then $f$ is associated with an homomorphism $g:\hz(1)\to G$. We have that $f$ is constant if and only if $g$ is trivial, or equivalently if and only if
    Using the image of the preferred $k$-rational section of $\cB_{k}\hz(1)$ corresponding to the trivial torsor, we may identify $\Pi_{T/k}$ with $\mc{B}_{k}\upi_{1}(T)$, and we reduce to prove the following: if $f:\hz(1)\to\pi_{1}(T_{\bar{k}})$ is a Galois equivariant homomorphism such that the composition $\H^{1}(k,\hz(1))\xar{f} \H^{1}(k,\pi_{1}(T_{\bar{k}})) \xar{res} \H^{1}(k_{\nu},\pi_{1}(T_{\bar{k}}))$ has finite image, then $f$ is trivial. 

    Choose $k'/k$ a finite extension such that $T_{k'}\simeq\G_{m}^{n}$, we get an identification $\pi_{1}(T_{\bar{k}})\simeq \hz^{n}(1)$ of Galois modules over $k'$. Let $\nu'$ be an extension of $\nu$ to $k'$, denote by $r$ the ramification index of $\nu'/\nu$. The valuation $\nu$ induces an homomorphism $\H^{1}(k_{\nu},\hz(1))=\projlim_{i}k_{\nu}^{*}/k_{\nu}^{*i}\to \hz$, and similarly for $k'_{\nu'}$; with an abuse of notation, we denote these homomorphisms by $\nu$, $\nu'$. 
    
    %The hypothesis implies that the composition 
    %\[\H^{1}(k,\hz(1)) \xar{f}  \H^{1}(k,\pi_{1}(T_{\bar{k}})) \xar{res} \H^{1}(k_{\nu},\pi_{1}(T_{\bar{k}})) \xar{res} \H^{1}(k'_{\nu'},\pi_{1}(T_{\bar{k}}))=\H^{1}(k'_{\nu'},\hz(1)^{n})\] has finite image.

    We have a commutative diagram
    \[\begin{tikzcd}
        \H^{1}(k,\hz(1))\rar["res"] \ar[dd,"res\circ \H^{1}(f)"]     &   \H^{1}(k_{\nu},\hz(1))\dar["res"]\rar["\nu"]    &   \hz\dar["r"]    \\
                                                &   \H^{1}(k'_{\nu'},\hz(1))\dar["\H^{1}(f)"]\rar["\nu'"]            &   \hz\dar["f"]         \\
        \H^{1}(k_{\nu},\pi_{1}(T_{\bar{k}}))\rar["res"] & \H^{1}(k'_{\nu'},\pi_{1}(T_{\bar{k}})\simeq\hz(1)^{n})\rar["\nu'"] & \pi_{1}(T_{\bar{k}})\simeq \hz^{n}
    \end{tikzcd}\]
    where the last identification $\pi_{1}(T_{\bar{k}})\simeq \hz^{n}$ ignores the Galois action. The hypothesis implies that the composition 
    \[f\circ r \circ \nu \circ res = \nu' \circ res \circ \H^{1}(f) : \H^{1}(k,\hz(1)) \to \hz^{n}\]
    has finite image, hence it is $0$ since $\hz$ has no torsion.  Notice that, since $\nu:k^{*}\to\Z$ is surjective, then the upper horizontal arrow $\nu \circ res: \H^{1}(k,\hz(1))=\projlim_{i}k^{*}/k^{*i} \to \hz$ is surjective as well, and hence $f\circ r =rf=0$. Since $r\neq 0$, we get that $f=0$, as desired.
    
%  It is enough to show that $f_{\ell}:\Z_{\ell}(1)\to \on{T}_{\ell}T=\projlim_{n}T[\ell^{n}]$ is trivial for every prime $\ell$. There exists a finite extension $h/k$ such that $T_{h}=\G_{m}^{n}$, let $f_{i}:\Z_{\ell}(1)\to\Z_{\ell}(1)$ be the $i$-th coordinate of $f_{\ell,h}$, it is enough to show that $f_{i}$ is trivial for every $i$.
		
%		Let $\mu$ be an extension of $\nu$ to $h$, $r$ the ramification index of $\mu/\nu$. Denote by $I_{\ell},J_{\ell}$ the projective limits $\projlim_{n}k^{*}/k^{*\ell^{n}}$, $\projlim_{n}h_{\mu}^{*}/h_{\mu}^{*\ell^{n}}$, we have $\H^{1}(k,\Z_{\ell})=I_{\ell}$, $\H^{1}(k_{\nu},\Z_{\ell})=J_{\ell}$ and the valuations $\nu,\mu$ induce a natural commutative diagram
%		\[\begin{tikzcd}
%			I_{\ell}\rar["\nu"]\dar	&	\Z_{\ell}\dar["r"]	\\
%			J_{\ell}\rar["\mu"]		&	\Z_{\ell}
%		\end{tikzcd}\]
		
%		If $\phi:J_{\ell}\to J_{\ell}$ is the homomorphism induced by $f_{i}$, the image of the composition
%		\[I_{\ell}\to J_{\ell}\xar{\phi}J_{\ell}\xar{\mu}\Z_{\ell}\]
%		is $rf_{i}(1)\Z_{\ell}$. The hypothesis implies that $I_{\ell}\to J_{\ell}\xar{\phi}J_{\ell}$ has finite image, hence $f_{i}(1)=0$ and $f_{i}$ is trivial.
	\end{proof}
\end{lemma}

\begin{lemma}\label{uniqueloop1}
	Let $k$ be a number field with a finite place $\nu$, $X$ a curve over $k$, $\gamma:\cB_{k}\hz(1)\to\Pi_{X/k}$ a morphism of gerbes. If the image of $\cB_{k}\hz(1)(k)$ in $\Pi_{X/k}(k_{\nu})$ has a finite number of isomorphism classes, then $\gamma$ is constant.
	\begin{proof}
		If by contradiction $\gamma$ is not constant, i.e. the associated homomorphism $\hz\to\pi_{1}(X_{\bar{k}})$ is not trivial, up to replacing $X$ with a finite étale cover we may assume that the composition $\cB_{k}\hz(1)\to\Pi_{X/k}\to\Pi_{X/k}^{\rm ab}$ is not constant, where $\Pi_{X/k}^{\rm ab}$ is the abelianized fundamental gerbe, equivalently we may assume that the induced homomorphism $\hz\to\pi_{1}(X_{\bar{k}})^{\rm ab}$ is not trivial. We may do so because the map of Galois sections associated with a finite étale cover has finite fibers, hence the hypothesis remains true after passing to the covering. 
		
		Let $J$ be the semi-abelian Jacobian of $X$, it is an extension of an abelian variety $A$ by a torus $T$ and we have a canonical identification $\Pi_{X/k}^{\rm ab}=\Pi_{J/k}$. We may use the images of preferred section of $\cB_{k}\hz(1)$ to give identifications $\Pi_{J/k}=\cB_{k}\on{T}J$, $\Pi_{A/k}=\cB_{k}\on{T}A$ where $\on{T}J=\projlim_{n}J[n],\on{T}A=\projlim_{n}A[n]$ are the Tate modules; the morphisms of gerbes then correspond to homomorphisms $\hz(1)\to\on{T}J$, $\on{T}J\to\on{T}A$.

		 %~ $\cB_{k}\hz(1)\to\Pi_{T/k}\to\Pi_{J/k}$ since $\Pi_{A/k}=\cB_{k}\on{T}A$.
		By weight reasons, the composition $\hz(1)\to\on{T}J\to\on{T}A$ is trivial, hence we have a factorization $\hz(1)\to\on{T}T\to\on{T}J$. We have a short exact sequence
		\[\projlim_{n}A[n](k_{\nu})\to\H^{1}(k_{\nu},\on{T}T)\to\H^{1}(k_{\nu},\on{T}J).\]
		By Mattuck's theorem \cite{mat55} $A(k_{\nu})$ has finite torsion, thus $\H^{1}(k_{\nu},\on{T}T)\to\H^{1}(k_{\nu},\on{T}J)$ is injective. The hypothesis implies that the image of $\H^{1}(k,\hz(1))$ in $\H^{1}(k_{\nu},\on{T}J)$ is finite, hence the image in $\H^{1}(k_{\nu},\on{T}T)$ is finite too. By \autoref{uniqueloop}, $\hz(1)\to\on{T}T$ is trivial, which is contradiction with the fact that $\cB_{k}\hz(1)\to\Pi_{X/k}^{\rm ab}$ is not constant.
		
		%~ $B\hz(1)(k)/{\sim}$ in $\Pi_{J/k}(k_{\nu})/{\sim}$ is finite. We have a short exact sequence
		%~ \[\projlim_{n}A[n](k_{\nu})\to\H^{1}(k_{\nu},\on{T}T)\to\H^{1}(k_{\nu},\on{T}J).\]
		%~ By Mattuck's theorem \cite{mat55} $A(k_{\nu})$ has finite torsion, thus $\Pi_{T/k}(k_{\nu})\to\Pi_{J/k}(k_{\nu})$ is injective on isomorphism classes and thus the image of $B\hz(1)(k)/{\sim}$ in $\Pi_{T/k}(k_{\nu})/{\sim}$ is finite, too. This implies that $B\hz(1)\to\Pi_{T/k}$ is constant by \autoref{uniqueloop}, hence we have a contradiction.
	\end{proof}
\end{lemma}

\section{Properties of $t$-birationally liftable Galois sections}\label{tprop}

%~ We call a morphism $X\to C$ a \emph{family of curves} if there exists a smooth, projective morphism $\bar{X}\to C$ with connected fibers of dimension $1$ and a divisor $D\s C$ finite étale over $C$ such that $X=\bar{X}\setminus D$. A family of curves is a geometric fibration, thus it has a relative étale fundamental gerbe $\Pi_{X/C}$.

%~ \subsection{Specializations of $t$-birational sections}

In order to work with $t$-b.l. Galois sections, we need to prove a number of facts about them.

\begin{lemma}
    Let $R$ be a DVR, $X\to \spec R$ a family of curves and $p\in X$ a closed point in the closed fiber. There exists a non-empty divisor $D\subset X$ finite étale over $R$ with $p\in D$.
	\begin{proof}
		%~ Let $K$, $k$ be the fraction and residue fields of $R$. First, let us show that it is enough to prove the statement in the case in which $p$ is $k$-rational. Let $k'$ be the residue field of $p$, since we are in characteristic $0$ we have $k'\simeq k[x]/\bar{q}(x)$ for some irreducible, separable polynomial $\bar{q}\in k[x]$, let $q\in R[x]$ any lifting of $\bar{q}$. By construction, $R'=R[x]/q(x)$ is finite étale over $R$ and has only one closed point, hence it is a DVR. If $D'\subset X_{R'}$ is finite étale over $R'$, it is finite étale over $R$, too, hence the image of $D'$ in $X$ is finite étale over $R$. 
		
		%~ Assume that $p$ is $k$-rational. By definition, $X\to\spec R$ has the form $\bar{X}\setminus D_{0}\to\spec R$ with $\bar{X}\to\spec R$ smooth, projective, with geometrically connected fibers of dimension $1$, and $D_{0}\subset\bar{X}$ a divisor finite étale over $R$. Choose a projective embedding $\bar{X}\subset\P^{n}_{R}$.
		
		%~ A generic hyperplane $H_{k}$ of $\P^{n}_{k}$ passing through $p$ has normal crossings with $\bar{X}_{k}$ and does not intersect $D_{0,k}$. Let $q\in R[t_{0},\dots,t_{n}]$ be a linear form whose image in $k[t_{0},\dots,t_{n}]$ is an equation for $H_{k}$, then $q$ defines a closed subscheme $H\subset\P^{n}_{R}$ which is a relative hyperplane. The fact that $H_{k}\cap D_{0,k}=\emptyset$ implies that $H\cap D_{0}=\emptyset$ since everything is proper over $R$, and the fact that $H_{k}$ has normal crossings with $\bar{X}_{k}$ implies that $D=H\cap X=H\cap\bar{X}$ is finite étale over $R$.

		%~ \newpage
		
		Let $K,k$ be the fraction and residue fields of $R$. If $p$ is $k$-rational, this essentially follows from Bertini's theorem by taking a generic hyperplane section passing through $p$. If $p$ is not $k$-rational, though, this is more subtle. We are going to use a Bertini-like argument to reduce to the case of $\P^{1}_{R}$, where we are able to make an explicit construction.
		
		By definition, $X\to\spec R$ has the form $\bar{X}\setminus D_{0}\to\spec R$ with $\bar{X}\to\spec R$ smooth, projective, with geometrically connected fibers of dimension $1$, and $D_{0}\subset\bar{X}$ a divisor finite étale over $R$. Choose a projective embedding $\bar{X}\subset\P^{n}_{R}$ and let $p_{1},\dots,p_{r}$ be the geometric points over $p$.
			
		A generic $(n-2)$-dimensional linear subspace $L_{k}$ of $\P^{n}_{k}$ has the following properties
		\begin{itemize}
			\item $L_{k}\cap \bar{X}_{k}=\emptyset$,
			\item for every $i$ the hyperplane $H_{i}\subset\P^{n}_{\bar{k}}$ spanned by $L_{k}$ and $p_{i}$ has transversal intersections with $\bar{X}_{\bar{k}}$,
			\item $\forall i:H_{i}\cap D_{0}=\emptyset$.
		\end{itemize}
		
		Let $\bar{q}_{0},\bar{q}_{1}\in k[t_{0},\dots,t_{n}]$ two $k$-linear forms which are equations for $L_{k}$, choose $q_{0},q_{1}\in R[t_{0},\dots,t_{n}]$ $R$-linear forms which lift $\bar{q}_{0},\bar{q}_{1}$ and let $L\subset\P^{n}_{R}$ be their vanishing locus, the fact that $L_{k}\cap \bar{X}_{k}=\emptyset$ implies that $L\cap \bar{X}=\emptyset$ since everything is proper over $R$.
		
		The graded homomorphism $R[x_{0},x_{1}]\to R[t_{0},\dots,t_{n}]$, $x_{i}\mapsto q_{i}$ defines a morphism 
		\[\P^{n}_{R}\setminus L=\proj R[t_{0},\dots,t_{n}]\setminus L\to\P^{1}_{R}=\proj R[x_{0},x_{1}].\]
		We have simply constructed the projection centered in $L$: we can do this even though we don't have a base field.
		
		Since $X\cap L=\emptyset$, we get a morphism $f:\bar{X}\to\P^{1}_{R}$ over $R$ with $\bar{X}_{k}\to\P^{1}_{k}$ surjective. The morphism $f$ is flat since it is finite and dominant, and both $\bar{X}$ and $\P^{1}_{R}$ are regular. Moreover, $f$ is unramified, and hence étale, at $\spec k(X_{k})\subset X$: notice that the diagram
		\[\begin{tikzcd}
			\bar{X}_{k}\rar\dar	&	\P^{1}_{k}\rar\dar	&	\spec k\dar	\\
			\bar{X}\rar			&	\P^{1}_{R}\rar		&	\spec R
		\end{tikzcd}\]
		is cartesian, hence $\spec k(\P^{1}_{k})\times_{\P^{1}_{R}}\bar{X}=\spec k(\bar{X}_{k})$ is reduced. By purity, it follows that the branch locus of $f$ is the closure of a finite number of closed points of $\P^{1}_{K}$. Let $K''/K$ be some finite extension and $H_{K''}\subset \P^{n}_{K''}$ a hyperplane containing $L_{K''}$ and corresponding to a branch point in $\P^{1}_{K}$; equivalently, $H_{K''}$ does not have transversal intersections with $X_{K''}$. If $H_{k''}\subset \P^{n}_{k''}$ is a specialization of $H_{K''}$ in some finite extension $k''$ of $k$, then $H_{k''}$ does not have transversal intersections with $X_{k''}$, it follows that $H_{k''}$ does not contain $p_{i}$ for every $i$. This implies that $H_{k''}$ corresponds to a point of $\P^{1}(k'')$ different from $f(p)$, hence $f(p)$ is not a branch point for $f$.
		
		Let $U_{0}\subset\P^{1}_{R}$ be the open subset where $f$ is étale, i.e. the complement of the image of the ramification locus in $\bar{X}$, and $U=U_{0}\setminus f(D_{0})$. Since $f(p)$ is not a branch point and $H_{i}\cap D_{0}=\emptyset$ for every $i$ we get that $f(p)\in U$. Furthermore, $f^{-1}(U)\to U$ is finite étale by construction. It is then enough to find a divisor of $U$ containing $f(p)$ and which is finite étale over $R$. 
		 %~ fact that $H_{i}$ has normal crossing with $\bar{X}_{\bar{k}}$ for every $i$ implies that $f(p)\in U_{0}$, while 
		
		%~ the fact that $H_{i}\cap D_{0}=\emptyset$ implies that $f(p)\in U$. Furthermore, $f^{-1}(U)\to U$ is finite étale by construction. It is then enough to find a divisor of $U$ containing $f(p)$ and which is finite étale over $R$.
		
		Since we are in characteristic $0$, then $k(p)\simeq k[x]/\bar{q}(x)$ for some irreducible, separable polynomial $\bar{q}\in k[x]$, let $q\in R[x]$ be any lifting of $\bar{q}$. By construction, $R'=R[x]/q(x)$ is finite étale over $R$ and has only one closed point, hence it is a DVR. Let $K'$ be the fraction field of $R'$, there are infinitely many points of $\P^{1}(K')$ which specialize to $p\in\P^{1}(k(p))$, hence we may choose one in $U(K')$; denote by $E\subset \P^{1}_{R}$ its closure, we have an extension $\spec R'\to E\subset\P^{1}_{R}$. By construction, $E$ is finite and flat over $R$ of degree equal to $[K':K]$. Furthermore, $E$ contains $p$, which satisfies $[k(p):k]=[K':K]$. Because of this, $E$ contains only $p$ and a generic point contained in $U$, hence $E\subset U$ and $\spec R'\to E$ is surjective. Since $R'$ is étale over $R$ we get that $E$ is étale over $R$ as well (we have actually proved that $\spec R'\to E$ is an isomorphism, though we don't need this).
  
%  The extension $\spec R'\to U$ is a closed embedding by degree reasons and defines the desired divisor.
		
		Notice that if we try to work with $R'$-sections of $\bar{X}$ we run into problems: it might be that $X(K')$ does not contain points which specialize to $p$, and if we try to enlarge $R'$ again we might get a semi-local ring with several closed points which we are not able to control simultaneously. Making a Bertini-like argument over $\bar{X}_{R'}$ runs into problems too: while we may find the desired divisor étale over $R'$, there is no guarantee that its image in $\bar{X}$ is étale over $R$.

	\end{proof}
\end{lemma}

\begin{corollary}\label{cut}
	Let $R$ be a DVR, $X\to \spec R$ a family of curves. There exists a direct system $(D_{i})_{i}$ of divisors finite étale over $R$ such that $X_{k}\setminus\bigcup_{i}D_{i}$ contains only the generic point of $X_{k}$.
\end{corollary}

\begin{lemma}\label{spebir}
	Specializations of b.l. sections are b.l.
	\begin{proof}
		Let $X\to \spec R$ be a family of curves with $R,K,k$ as above and $c\in\spec R$ the closed point. If $z\in\Pi_{X_{K}/K}(K)$ is a generic, b.l. Galois section, we have a specializing loop $\gamma_{z}:S^{1}_{c}\to\Pi_{X_{k}/k}$. 
		%~ We want to prove that the sections of $\Pi_{X_{p}/k(p)}(k(p))$ in the essential image of $\gamma_{z}(k(p))$ are birational.
		
		Let $D_{i}\s X$ be a direct system of divisors as in \autoref{cut}. We have that $X_{i}=X\setminus D_{i}$ is a family of curves for every $i$, write $X_{\infty}=\projlim_{i}X_{i}$ and $\Pi_{X_{\infty}/R}=\projlim_{i}\Pi_{X_{i}/R}$. The closed fiber of $X_{\infty}$ is naturally isomorphic to $\spec k(X_{k})$.
		
		Since $z$ is b.l. and $\spec K(X)\s X_{\infty,K}$, there exists a lifting $z'\in\Pi_{X_{\infty,K}/K}(K)$ of $z$ to $X_{\infty,K}$. By a limit argument, the specializing loops of $z'$ in $\Pi_{X_{i}/R}$ induce a specializing loop $\gamma_{z'}:S^{1}_{c}\to \Pi_{X_{\infty,k}/k}=\Pi_{k(X_{k})/k}$, and the composition of $\gamma_{z'}$ with $\Pi_{k(X_{k})/k}\to \Pi_{X_{k}/k}$ is isomorphic to $\gamma_{z}$. Hence, every specialization of $z$ lifts to $\spec k(X_{k})$.
	\end{proof}
\end{lemma}

\begin{corollary}\label{tbb}
	 $t$-b.l. sections are b.l. 
	 \begin{proof}
		Let $X$ be a curve over a field $k$ and $s\in\Pi_{X/k}(k)$ a $t$-b.l. Galois section. Write $R=k[t]_{(t)}$, we have that $s_{k(t)}$ is a generic section of $\Pi_{X_{R}/R}$ and $s$ is a specialization of $s_{k(t)}$. By hypothesis $s_{k(t)}$ is b.l., hence its specialization $s$ is b.l. by \autoref{spebir}. 
	 \end{proof}
\end{corollary}

\begin{lemma}\label{spetbir}
	Specializations of $t$-b.l. sections are $t$-b.l.
	\begin{proof}
		Let $X\to \spec R$ be as above, $z\in\Pi_{X_{K}/K}(K)$ a generic, $t$-b.l. Galois section, $s\in\Pi_{X_{k}/k}(k)$ a specialization. Let $R'$ be the local ring of the generic point of the divisor $\A^{1}_{k}\s\A^{1}_{R}$. We have that $R'$ is a DVR with fraction field $K'=K(t)$ and residue field equal $k'=k(t)$. Consider the family of curves $X_{R'}\to \spec R'$, then $z_{K(t)}=z_{K'}\in \Pi_{X_{K'}/K'}(K')$ is by hypothesis a b.l. Galois section. We have that $s_{k(t)}=s_{k'}$ is a specialization of $z_{K'}$ and thus it is b.l. by \autoref{spebir}.
	\end{proof}
\end{lemma}

%~ \subsection{Lifting $t$-birational sections}

\begin{lemma}\label{morb}
	Let $f:Y\to X$ be a dominant morphism of curves over a field $k$ and $s\in\Pi_{Y/k}(k)$ a Galois section. If $s$ is b.l. then $f(s)\in\Pi_{X/k}(k)$ is b.l. If $f$ is finite étale, the converse holds.
	\begin{proof}
		The first statement is obvious. If $f$ is finite étale, the second statement follows from the fact that $\gal(k(Y))=\pi_{1}(Y)\times_{\pi_{1}(X)}\gal(k(X))$.
	\end{proof}
\end{lemma}

\begin{corollary}\label{mortb}
	Let $f:Y\to X$ be a dominant morphism of curves over a field $k$ and $s\in\Pi_{Y/k}(k)$ a Galois section. If $s$ is $t$-b.l. then $f(s)\in\Pi_{X/k}(k)$ is $t$-b.l. If $f$ is finite étale, the converse holds.\qed
\end{corollary}

We say that a curve $X$ is \emph{parabolic} if it has genus $0$ and the degree of $\bar{X}\setminus X$ is at most $2$. If $X$ is proper of genus $0$, then there is a unique Galois section, and it is geometric if and only if $X=\P^{1}$. If $X$ is parabolic affine, every Galois section of $X$ is cuspidal (the only non-trivial case is settled by \cite[Theorem A]{sch15}). Because of this, we can always assume that $X$ is non-parabolic.

If $X$ is a curve over a number field $k$ and $s$ is a b.l. section, by Koenigsmann's results \cite{koe05}, \cite[Proposition 1]{sti15} for every finite place $\nu$ of $k$ there exists a local point $x_{\nu}\in\bar{X}(k_{\nu})$ such that $s_{k_{\nu}}$ is associated with $x_{\nu}$, i.e. if $x_{\nu}\in X$ then $s_{\nu}$ is the geometric section associated with $x_{\nu}$, otherwise $s_{\nu}$ is one of the cuspidal sections of the packet associated with $x_{\nu}\in \bar{X}\setminus X$.

If $X$ is non-parabolic, the local point $x_{\nu}$ associated with $s$ is unique. To check this, we can pass to an étale neighbourhood $Y$ of genus $\ge 1$ and use the injectivity of the section map \cite[Proposition 75]{sti13} on $\bar{Y}_{k_{\nu}}$.

%~ \begin{lemma}
	%~ Let $k$ be a finite extension of $\Q_{p}$, $X$ a curve over $k$ with good reduction. Suppose we have a morphism $h:B\hz(1)\to\Pi_{X/k}$ such that there are only finitely many isomorphism classes in the essential image of $B\hz(1)(k)\to\Pi_{X/k}(k)$. Then $h$ is constant.
	%~ \begin{proof}
			%~ If $h$ is not constant, up passing to an étale cover of $X$ we may assume that the composition $h^{\rm ab}:B\hz(1)\to\Pi_{X/k}\to\Pi_{X/k}^{\rm ab}$ is not constant, where $\Pi_{X/k}^{\rm ab}$ is the abelianization.
	%~ \end{proof}
%~ \end{lemma}

\begin{lemma}\label{nflift}
	Let $X$ be a curve over a number field $k$, $U\s X$ a non-empty open subset and $s\in\Pi_{X/k}(k)$ a $t$-b.l. Galois section. There exists a $t$-b.l. section of $U$ which lifts $s$.
	\begin{proof}
		We may assume that $X$ is non-parabolic and that $s$ is neither geometric nor cuspidal, otherwise this is trivial. 
		
		By hypothesis, $s_{k(t)}$ lifts to a b.l. Galois section $r_{0}$ of $U_{k(t)}$, we are going to show that $r_{0}$ is the base change to $k(t)$ of some Galois section of $U$. Identify $k(t)$ with the fraction field of $\A^{1}$ and choose $c\in \A^{1}$ a closed point with residue field $k'$.
  
        For every finite place $\nu$ of $k$ there is a unique local point $x_{\nu}\in\bar{X}(k_{\nu})$ associated with $s$. Since $s$ lifts to a b.l. section of $U$ by \autoref{tbb} and it is not geometric nor cuspidal, the set of places $\nu$ such that $x_{\nu}\in U$ has Dirichlet density $1$ by \cite[Theorem B]{sti15}, hence we may choose $\nu$ with $x_{\nu}\in U$ and $k'\subset k_{\nu}$. If $r_{0,c}$ is a specialization of $r_{0}$ at $c$, by construction it is b.l. and lifts $s$, hence $r_{0,c,k_{\nu}}$ is geometric associated with $x_{\nu}$. By \autoref{uniqueloop1}, the specializing loop $S^{1}_{c}\to\Pi_{U_{k'}/k'}$ of $r_{0}$ at $c$ is constant and hence $r_{0}$ extends to a section $\spec\O_{c}\to\Pi_{U\times\A^{1}/\A^{1}}$ by \autoref{prorootval}.
		
		Since we can do this for every closed point $c\in\A^{1}$ and we may write $\Pi_{U\times\A^{1}/\A^{1}}$ as a product $\Pi_{U/k}\times_{k}\A^{1}$, by \cite[Corollary A.3]{fed} we get that $r_{0}$ extends to a section $\tilde{r}_{0}:\A^{1}\to\Pi_{U\times\A^{1}/\A^{1}}$. Since $\Pi_{\A^{1}/k}=\spec k$, the composition $\A^{1}\to\Pi_{U\times\A^{1}/\A^{1}}\to\Pi_{U/k}$ factorizes through a section $r:\spec k\to\Pi_{U/k}$ and hence $\tilde{r}_{0}=r_{\A^{1}}$, $r_{0}=r_{k(t)}$. This implies that $r$ is a $t$-b.l. lifting of $s$.
	\end{proof}
\end{lemma}

\section{The main argument}

This section consists of a unique statement. The main theorems will follow easily from it.

\begin{proposition}\label{simple}
	Let $k$ be a number field, $\nu$ a finite place, write $U=\A^{1}\setminus\{1,2\}$. Let $s\in\Pi_{U/k}(k)$ be a b.l. Galois section and $x_{\nu}\in\P^{1}(k_{\nu})$ the unique associated local point. Let $\delta:\spec k(t)\to\A^{1}$ be the ``diagonal''	 point, i.e. the generic one, and assume that $s_{k(t)}$ lifts to a section of $\spec k(\A^{1})\otimes_{k} k(t)\setminus\{\delta\}$ (e.g. if $s$ is $t$-b.l.). Then $x_{\nu}$ is $k$-rational.
	\begin{proof}
		We may assume that $x_{\nu}\in U$ since otherwise it is clearly $k$-rational. Let $r_{0}$ be the lifting of $s_{k(t)}$ to $\spec k(\A^{1})\otimes_{k} k(t)\setminus\{\delta\}$, and $r$ the image of $r_{0}$ in $U_{k(t)}\setminus \{\delta\}$. We divide the proof in three steps.
		
		\textbf{Step 1. Specializations of $r$.} We are going to prove that, for every $k$-rational point $c\in U(k)$ different from $x_{\nu}$, the specializing loop of $r$ is constant: the idea is to prove that for any specialization $r_{c}$ of $r$ the base change of $r_{c}$ to $k_{\nu}$ is the geometric section associated with $x_{\nu}\in U\setminus\{c\}(k_{\nu})$ regardless of the choice of $r_{c}$, then we apply \autoref{uniqueloop1}.
		
		Write $\Delta\s U\times U$ for the diagonal, consider the family of curves $U\times U\setminus \Delta\to U$ and fix a $k$-rational point $c\in U(k)$ with $c\neq x_{\nu}$, we have $(U\times U\setminus \Delta)_{c}=U\setminus\{c\}$. Consider the scheme $W=\spec\O_{c}\otimes_{k}\O_{c}\setminus\Delta$ over $\spec \O_{c}$, it contains $\spec k(\A^{1})\otimes_{k} k(t)\setminus\{\delta\}$. The morphism $W\to\spec\O_{c}$ is a projective limit of families of curves over $\O_{c}$ (it can be obtained from $U_{\O_{c}}\to\spec \O_{c}$ by removing divisors étale over $\O_{c}$), its generic fiber is $\spec \O_{c}\otimes_{k} k(t)\setminus\{\delta\}$ while the special fiber is $\spec\O_{c}\setminus\{c\}=\spec k(t)$. 
		
		Choose $r_{c}\in\Pi_{U\setminus\{c\}/k}(k)$ any specialization of $r$. Let $r_{1}$ be the image of $r_{0}$ in $\spec \O_{c}\otimes_{k} k(t)\setminus\{\delta\}$. The specializing loop $S^{1}_{c}\to\Pi_{U\setminus\{c\}/k}$ of $r$ factorizes through the specializing loop $S^{1}_{c}\to\Pi_{k(t)/k}$ of $r_{1}$, hence $r_{c}$ is b.l. Furthermore, the image of $r_{c}$ in $U$ is a specialization of $s_{k(t)}$, i.e. it is $s$.
		
		Since $x_{\nu}\in U\setminus\{c\}$, $r_{c}$ is b.l. and maps to $s$, we have that the base change $r_{c,k_{\nu}}$ is \emph{the} geometric section associated to $x_{\nu}$, regardless of the choice of $r_{c}$. By \autoref{uniqueloop1}, the specializing loop of $r$ at $c$ is constant.
	
		%%%%%
	
		\textbf{Step 2. Change of coordinates.} The map $y\mapsto t-y$ defines an automorphism $\phi:\A^{1}_{k(t)}\to\A^{1}_{k(t)}$ with $\phi(\delta)=0$. Write $V=\A^{1}_{k(t)}\setminus\{\phi(1),\phi(2)\}$, we have that $\phi$ restricts to an isomorphism $\phi:U_{k(t)}\to V$. We thus get a Galois section $\phi(s_{k(t)})\in\Pi_{V}(k(t))$ with a lifting $\phi(r)\in\Pi_{V\setminus \{0\}}(k(t))$. Since $V\setminus\{0\}\s\G_{m}$, then $\phi(r)$ induces a section
		\[z\in\Pi_{\G_{m}/k(t)}(k(t))=\cB\hz(1)(k(t))=\projlim_{n} k(t)^*/k(t)^{*n}=\hat{k(t)^*},\]
		\[z=\lambda\cdot\prod_{c}q_{c}^{e_{c}}\]
		where $\lambda\in\hat{k^*}$, $c$ varies among closed points of $\A^{1}$, $q_{c}\in k[t]$ is the monic, irreducible polynomial associated with $c$ and $e_{c}\in\hz$.\footnote{Notice that, while there is an embedding $\hat{k(t)^{*}}\subset \hat{k^{*}}\cdot\prod_{c}q_{c}^{\hz}$, this is not a bijection, since for every $n$ all but finitely many exponents of an element of $\hat{k(t)^{*}}$ are multiples of $n$.} For every $k$-rational point $c\neq x_{\nu},1,2$, the fact that the specializing loop of $r$ at $c$ is constant implies that the same holds for $\phi(r)$ and hence for $z$, hence $e_{c}=0$ by \autoref{gmloop}. Since we can repeat step 1 after base changing to any finite extension of $k$, we get that $e_{c}=0$ for every \emph{closed} point $c\neq x_{\nu},1,2$.
		
		If $x_{\nu}$ is algebraic over $k$, let $q$ be its minimal polynomial and $e$ the exponent of $q$ as a factor of $z$, otherwise $q=1$ and $e=0$. We may thus write
		\[z=\lambda\cdot q^{e}\cdot(t-1)^{e_{1}}\cdot(t-2)^{e_{2}}.\]
		We are going to prove that $q$ has degree $1$ and thus $x_{\nu}$ is $k$-rational.
		
		%%%%%%
		
		\textbf{Step 3. Specializations of $z$.} Fix $K$ a finite extension of $k$ which splits $q$ completely (if $q=1$ choose $K=k$) and let $\eta$ be an extension of $\nu$. In the rest of the proof, we are going to consider many $K$-rational points $c\in K$ (or $c_{i}\in K$). We will always tacitly assume that $c\neq 1,2$ and $q(c)\neq 0$.
		
		If we repeat steps 1 and 2 after base changing to $K$, for every $c$ the base change to $K_{\eta}$ of the unique specialization of $\phi(r)$ at $c$ is $\phi(x_{\eta})$, hence we obtain the equation
		\[\lambda\cdot q(c)^{e}\cdot(c-1)^{e_{1}}\cdot(c-2)^{e_{2}}=c-x_{\eta}\in\hat{K_{\eta}}^{*},\]
		and by applying $\eta:\hat{K_{\eta}}^{*}\to\hz$ we get
		\[\eta(\lambda)+e\eta(q(c))-\eta(c-x_{\eta})+e_{1}\eta(c-1)+e_{2}\eta(c-2)=0\in\hz.\]
		
		Observe that if $p\in K[t]$ is a polynomial, $c\in K$ an element with $p(c)\neq 0$ and $c_{i}\neq c$ is a sequence which tends to $c$ in the $\eta$-adic topology, then $\eta(p(c_{i}))=\eta(p(c))$ is constant for $i>>0$ large enough while $\eta(c_{i}-c)$ is \emph{not} constant.
		
		Choose a sequence of $K$-rational points $c_{i}$ which tends to $1$. Since $1$ is not a root of the polynomials $q$, $t-2$, $t-x_{\eta}$ we see that for $i>>0$ all terms except $e_{1}\eta(c_{i}-1)$ in the equation above are constant. It follows that $e_{1}\eta(c_{i}-1)$ is constant, too. Since $c_{i}\xar{i} 1$ and $c_{i}\neq 1$, then $\eta(c_{i}-1)$ is not constant, this implies that $e_{1}=0$. With the same argument we see that $e_{2}=0$, hence
		\[\eta(\lambda)+e\eta(q(c))-\eta(c-x_{\eta})=0\in\hz.\]
		
		If $x_{\nu}$ is transcendental over $k$ and thus $q=1$, $e=0$, then $\eta(c-x_{\eta})=\eta(\lambda)$ does not depend on $c$, which is absurd, hence $x_{\nu}$ is algebraic over $k$. Since $K$ splits $q$, we may write $q(t)=(t-x_{\eta})\cdot\prod_{j}(t-y_{j})$,
		\[\eta(\lambda)+(e-1)\eta(c-x_{\eta})+e\sum_{j}\eta(c-y_{j})=0\in\hz.\]
		Since we are in characteristic $0$ and $q$ is irreducible over $k$, then $x_{\nu}\neq y_{j}$ for every $j$ and $y_{j}\neq y_{j'}$ for $j\neq j'$. Using a sequence $c_{i}\xar{i}x_{\eta}$ and the same argument as above we see that $e=1$ and hence 
		\[\eta(\lambda)+\sum_{j}\eta(c-y_{j})=0\in\hz.\]
			
		If by contradiction $x_{\nu}$ is not $k$-rational and thus $\deg q\ge 2$, using a sequence $c_{i}\xar{i}y_{j}$ we see that $\eta(\lambda)+\sum_{j}\eta(c-y_{j})=0$ is not constant for $i>>0$, which is absurd.
	\end{proof}
\end{proposition}

\section{Reduction to number fields}

In \cite{saty21} M. Sa\"idi and M. Tyler reduced the birational section conjecture to number fields. We are going to do this for the $t$-birational version, too. A reader only interested in number fields may safely skip this section.

Over number fields, we know that we can lift $t$-b.l. sections to open subsets thanks to \autoref{nflift}. In order to prove it, we use a theorem of J. Stix \cite[Theorem B]{sti13} which is only available for number fields. Lifting to open subsets is crucial for the reduction to number fields, and we can only do it for number fields: this is a problem. In order to overcome it, we define \emph{quasi-$t$-b.l.} sections.

\begin{definition}
	Let $X$ be a geometrically connected curve over a field $k$. A Galois section $s$ of $X$ is quasi-$t$-b.l. if there is an open subset $U\s X$, a lifting $s'$ of $s$ to $U$ and a dominant map $f:U\to Y$ where $Y$ is a non-parabolic curve such that $s'$ is b.l. and $f(s')$ is $t$-b.l.
\end{definition}

Since $t$-b.l. sections are b.l. by \autoref{tbb}, it is immediate to check that $t$-b.l. sections are quasi-$t$-b.l. If $s$ is a quasi-$t$-b.l. section of $X$ and $U\s X$ is an open subset, then it is obvious that we can lift $s$ to a quasi-$t$-b.l. section of $U$.

\begin{lemma}\label{spequasi}
	Let $C$ be a curve with fraction field $K$, $X\to C$ a family of curves and $s\in\Pi_{X_{K}/K}(K)$ a generic section which is quasi-$t$-b.l. There exists a non-empty open subset $U\s C$ such that the specializations of $X|_{U}\to U$ are quasi-$t$-b.l.
	\begin{proof}
		Up to shrinking both $C$ and $X$, we may assume that there exists a family of non-parabolic curves $Y\to C$ and a quasi-finite morphism $f:X\to Y$ over $C$ such that $f(s)$ is $t$-b.l. The specializations of $f(s)$ are $t$-b.l. by \autoref{spetbir} and the specializations of $s$ are b.l. by \autoref{spebir}, it follows that they are quasi-$t$-b.l., too.
	\end{proof}
\end{lemma}

\begin{lemma}\label{quasinf}
	Let $k$ be a number field, and assume that $t$-b.l. sections over $k$ are geometric or cuspidal. Then the same holds for quasi-$t$-b.l. sections.
	\begin{proof}
		Let $X$ be a curve over $k$ and $s\in\Pi_{X/k}(k)$ a quasi-$t$-b.l. section. We may assume that there exists a non-parabolic curve $Y$ and a dominant morphism $f:X\to Y$ such that $f(s)$ is $t$-b.l. Since $s$ is b.l., for every finite place $\nu$ of $k$ there exists a unique associated point $x_{\nu}\in \bar{X}(k_{\nu})$. Since $f(s)$ is $t$-b.l., there exists a unique $k$-rational point $y\in \bar{Y}(k)$ associated to $f(s)$, in particular $f(x_{\nu})=y$ for every $\nu$. Any b.l. lifting of $s$ to $X\setminus f^{-1}(y)$ is cuspidal thanks to \cite[Theorem B]{sti15}, hence $s$ is geometric or cuspidal.
	\end{proof}
\end{lemma}

Let us recall a famous result by A. Tamagawa.

\begin{proposition}[Tamagawa]\label{tama}
    Let $X$ be a hyperbolic curve over a field $k$ finitely generated over $\Q$ and let $s\in\Pi_{X/k}(k)$ be a Galois section. If $s$ is not geometric nor cuspidal, there exists a curve $Y$ of genus $\ge 2$ with $\bar{Y}(k)=\emptyset$, a finite étale morphism $Y\to X$ and a lifting $r\in\Pi_{Y/k}(k)$ of $s$.
\end{proposition}

\begin{proof}
    This is essentially \cite[Proposition 2.8 (iv)]{tam97}.
\end{proof}

\begin{proposition}\label{nffg}
	Assume that $t$-b.l. sections of curves defined over number fields are geometric or cuspidal. Then quasi-$t$-b.l. sections of curves over fields of finite type over $\Q$ are geometric or cuspidal.
	\begin{proof}
		We prove this by induction on the transcendence degree $n$ of the base field $k$ over $\Q$. The case $n=0$ is \autoref{quasinf}. Assume $n \ge 1$ and that the statement is proved for fields of transcendence degree $n-1$, let $h\s k$ be algebraically closed in $k$ and of transcendence degree $n-1$ over $\Q$. Let $X$ be a curve and $s\in\Pi_{X/k}$ a quasi-$t$-b.l. section, up to passing to an open subset we may assume that $X$ is hyperbolic. Assume by contradiction that $s$ is not geometric nor cuspidal, thanks to Tamagawa's argument (\autoref{tama}) and \autoref{morb} up to passing to an étale neighbourhood we may assume that $\bar{X}(k)=\emptyset$. Clearly the image in $s$ in $\bar{X}$ is not geometric, so we may replace $X$ with $\bar{X}$ and assume that $X$ is projective.
		
		Choose $E$ any elliptic curve over $h$, we may find a finite, possibly ramified cover $Y\to X$ with a finite morphism $Y\to E_{k}$. Let $U\s X,V\s Y$ be open subsets such that $V\to U$ is finite étale. Since $s$ is quasi-$t$-b.l., by applying \autoref{morb} to $V\to U$ we may find a finite extension $k'/k$ and a quasi-$t$-b.l. lifting $r\in\Pi_{Y/k}(k')$ of $s_{k'}$. Thanks to \cite[Lemma 6.5, Theorem 7.2]{anab}, it is enough to prove that $r$, and thus $s_{k'}$, is geometric. We may thus replace $k$, $X$, $s$ with $k'$, $Y_{k'}$, $r$ and assume that there is a finite morphism $X\to E_{k}$.
		
		Let $C$ be an affine curve over $h$ whose function field is $k$, up to shrinking $C$ we may extend $X$ to a family of smooth projective curves $\tilde{X}\to C$ with a finite morphism $\tilde{X}\to E\times C$ over $C$. The generic section $s$ extends to a global section $\tilde{s}:C\to\Pi_{\tilde{X}/C}$ thanks to \cite[Corollary 3.4]{scfg}. The specializations of $\tilde{s}$ are geometric thanks to \autoref{spequasi} plus inductive hypothesis. It follows that $s$ is geometric thanks to \cite[Definition 4.1, Corollary 4.9]{scfg}.
	\end{proof}
\end{proposition}

\section{Proof of the main theorems}

Theorems A, B and C follow rather easily from \autoref{simple}.

\begin{TA}
	Let $X$ be a smooth curve over a field $k$ finitely generated over $\Q$. A Galois section of $X$ is geometric or cuspidal if and only if it is $t$-birationally liftable.
	\begin{proof}
		Thanks to \autoref{nffg}, we may assume that $k$ is a number field. Let $s\in\Pi_{X/k}(k)$ be a $t$-b.l. section, assume by contradiction that $s$ is neither geometric nor cuspidal. By \autoref{nflift} we may assume that $X$ is hyperbolic. By Tamagawa's argument (\autoref{tama}) and \autoref{mortb}, up to passing to an étale neighbourhood of $s$ we may assume that $X$ has genus at least $2$ and $\bar{X}(k)=\emptyset$. Fix $\nu$ any finite place of $k$ and let $x_{\nu}\in \bar{X}(k_{\nu})$ the local point associated with $s$.
		
		Choose a projective embedding $j:\bar{X}\s\P^{n}$ such that $j(x_{\nu})\in\A^{n}\s\P^{n}$ and let $U=\A^{n}\cap X$, we have that $s$ lifts to a $t$-b.l. section $s'$ of $U$ thanks to \autoref{nflift}. Since $\bar{X}(k)=\emptyset$ then $j(x_{\nu})\in\A^{n}$ is not $k$-rational, in particular there exists one coordinate $c:\A^{n}\to\A^{1}$ such that $c(j(x_{\nu}))$ is not $k$-rational. Up to shrinking $U$ furthermore we may assume that $c(j(U))\s\A^{1}\setminus\{1,2\}$. Then $c(j(s'))\in\Pi_{\A^{1}_{k}\setminus\{1,2\}/k}(k)$ is a $t$-b.l. section whose base change to $k_{\nu}$ is associated with a non-$k$-rational local point, which is in contradiction with \autoref{simple}.
	\end{proof}
\end{TA}

% First, let us show that the cuspidalization conjecture is equivalent to the statement that every Galois section is birationally liftable.

\begin{lemma}\label{cuspbl}
    Let $X$ be a curve over a countable field $k$. The following are equivalent.
    \begin{itemize}
        \item For every pair of non-empty open subset $V\subset U\subset X$, the map $\mc{S}_{V/k}\to\mc{S}_{U/k}$ is surjective.
        \item For every non-empty open subset $U\subset X$, every Galois section of $U$ is birationally liftable.
    \end{itemize}
\end{lemma}

\begin{proof}
    The second condition clearly implies the first. Assume that the first holds and let $X'\subset X$ be a non-empty open subset. Since $k$ is countable, then $X'$ has a countable number of open subsets. Since $\mc{S}_{k(X')/k}\simeq\projlim_{U\subset X'}\mc{S}_{U/k}$ \cite[Lemma 259]{sti13}, the projective system $(\mc{S}_{U/k})_{U}$ is countable and the transition maps are surjective, then $\mc{S}_{k(X')/k}\to\mc{S}_{X'/k}$ is surjective as well, i.e. every Galois section of $X'$ is birationally liftable.
\end{proof}

\begin{corollary}\label{cuspbl2}
    The following are equivalent.
    \begin{itemize}
        \item The cuspidalization conjecture holds.
        \item For every smooth, geometrically connected hyperbolic curve $X$ over a field $k$ finitely generated over $\Q$, every Galois section of $X$ is birationally liftable.
    \end{itemize}
\end{corollary}

\begin{TB}\label{cuspeq}
	Let $X$ be a hyperbolic curve over a field $k$ finitely generated over $\Q$. The following are equivalent.
	\begin{itemize}
		\item For every finitely generated extension $K/k$, every Galois section of $X_{K}$ is either geometric or cuspidal.
		\item For every finitely generated extension $K/k$, every Galois section of $X_{K}$ is birationally liftable.
	\end{itemize}
        As a consequence, the section conjecture is equivalent to the cuspidalization conjecture.
        
	\begin{proof}
		Follows directly from \hyperlink{TA}{Theorem A} and \autoref{cuspbl2}.
	\end{proof}
\end{TB}

\begin{lemma}\label{hilbsp}
	Let $k$ be a Hilbertian field, $V$ an open subset of $\A^{1}$, $X\to V$ a family of curves, $s_{1},s_{2}\in\Pi_{X/V}(V)$ two sections. If $s_{1},s_{2}$ have isomorphic specializations at $k$-rational points of $V$, then $s_{1}\simeq s_{2}$.
	\begin{proof}
		For the convenience of the reader, we give proofs in both the language of fundamental gerbes and fundamental groups.
		 
		\subsubsection*{With fundamental gerbes} Write $\Pi_{X/V}=\projlim\Phi_{i}$ as in \autoref{limit}, let $W_{i}$ be the fibered product $V\times_{\Phi_{i}}V$ with respect to the two morphisms $s_{1},s_{2}:V\to\Pi_{X/V}\to\Phi_{i}$, we have that $\uisom(s_{1},s_{2})=V\times_{\Pi_{X/V}}V=\projlim_{i}W_{i}$.
		
		There is a natural map $W_{i}=V\times_{\Phi_{i}}V\to V\times_{V}V=V$ which coincides with projection on both coordinates, it is a finite étale cover since $\Phi_{i}$ is a proper étale gerbe over $V$. By hypothesis, $\uisom(s_{1},s_{2})(k)\to V(k)$ is surjective, this implies that $W_{i}(k)\to V(k)$ is surjective too. Since $k$ is Hilbertian and $W_{i}\to V$ is finite étale, there is a section $V\to W_{i}$. Let $W_{i}(V)$ be the set of sections, it is finite and non-empty. It follows that $\uisom(s_{1},s_{2})(V)=\projlim_{i}W_{i}(V)$ is non-empty since a projective limit of finite, non-empty sets is non-empty. This gives the desired global isomorphism $s_{1}\simeq s_{2}$.
		
		\subsubsection*{With fundamental groups} For every $v\in V(k)$, choose $r_{v}:\gal(\bar{k}/k)\to\pi_{1}(V)$ a representative of the geometric Galois section associated with $V$. If $S\subset V(k)$ is not thin in the sense of Serre, the images of the sections $r_{v}$ for $v\in S$ generate $\pi_{1}(V)$ topologically: if $H\subset\pi_{1}(V)$ is an open subgroup containing all of them, the corresponding finite étale covering $C\to V$ is connected and the image of $C(k)\to V(k)$ contains $S$, hence $C=V$ since $S$ is not thin.
		
		Let $F$ be a geometric fiber of the family, we have that $\pi_{1}(X)$ is an extension of $\pi_{1}(V)$ by $\pi_{1}(F)$. We may write $\pi_{1}(X)$ as a projective limit $\projlim_{i}G_{i}$ of pro-finite groups such that $G_{i}$ is an extension of $\pi_{1}(V)$ by a finite group $H_{i}$ and $\pi_{1}(F)=\projlim_{i}H_{i}$; this can be done by writing first $\pi_{1}(X)=\projlim_{i}G_{0,i}=\pi_{1}(X)/N_{i}$ as a projective limit of finite groups and then defining $G_{i}=G_{0,i}\times_{\pi_{1}(V)/N_{i}}\pi_{1}(V)$. Choose representatives $t_{1},t_{2}:\pi_{1}(V)\to\pi_{1}(X)$ of $s_{1},s_{2}$, and let $t_{1,i},t_{2,i}:\pi_{1}(V)\to \pi_{1}(X)\to G_{i}$ be the compositions.
		
		For every $v\in V(k)$ and every $i$, by hypothesis there exists $h\in H_{i}$ such that $h^{-1}(t_{1,i}\circ r_{v})h=t_{2,i}\circ r_{v}$. Since $H_{i}$ is finite there exists an $h\in H_{i}$ such that the set $S_{h}\subset V(k)$ of points $v$ for which $h^{-1}(t_{1,i}\circ r_{v})h=t_{2,i}\circ r_{v}$ is not thin. Since the homomorphisms $r_{v}$ for $v\in S_{h}$ generate $\pi_{1}(V)$, this implies that $h^{-1}t_{1,i}h=t_{2,i}$.	Since a projective limit of finite, non-empty subsets is non-empty, we get an element $g\in\pi_{1}(F)=\projlim_{i}H_{i}$ such that $g^{-1}t_{1}g=t_{2}$, i.e. the $\pi_{1}$ sections $s_{1},s_{2}$ are equal.
	\end{proof}
\end{lemma}

%~ \begin{lemma}\label{hilbgen}
	%~ Let $k$ be a Hilbertian field, $V\s\P^{1}$ an open subset, $s_{v}:\gal(k)\to\pi_{1}(V)$ a choice of a section associated with $v\in V(k)$ for every $v$. The images of the sections $s_{v}$ generate $\pi_{1}(V)$ topologically.
	%~ \begin{proof}
		%~ Let $H\s\pi_{1}(V)$ be an open subgroup containing the image of $s_{v}$ for every $v$, it is associated with a finite étale morphism $C\to V$ with $C$ connected. Since $s_{v}$ maps to $H$, then there exists a rational point of $C$ over $v$. Since this is true for every $v\in V(k)$ and $k$ is Hilbertian, it follows that $C=V$ and $H=\pi_{1}(V)$.
	%~ \end{proof}
%~ \end{lemma}

Let us prove that the two forms of the birational section conjecture given in the introduction are equivalent.

\begin{lemma}\label{bireq}
    The following are equivalent.
    \begin{itemize}
        \item The birational section conjecture holds.
        \item For every curve $X$ over a field $k$ finitely generated over $\Q$, every b.l. Galois section of $X$ is either geometric or cuspidal.
    \end{itemize}
\end{lemma}

\begin{proof}
    The first condition clearly implies the second. Assume that the second holds, let $X$ be a curve over a field $k$ finitely generated over $\Q$ and let $s\in\mc{S}_{k(X)/k}$ be a birational Galois section, we want to show that $s$ is cuspidal. Up to replacing $X$ with an open subset, we can assume that $X$ is hyperbolic.
    
    For every open subset $U\subset X$, denote by $\iota_{U}(s)$ the image of $s$ in $\mc{S}_{U/k}$, by hypothesis it is cuspidal. By \cite[Theorem 17 (2)]{sti12}, there exists a unique rational point $x\in\bar{X}(k)$ associated with $\iota_{X}(s)$, where $\bar{X}$ is the smooth completion of $X$. By uniqueness, $x$ is associated with $\iota_{U}(s)$ for every open subset $U$ as well.
    
    For every open subset $U\subset X$, denote by $\mc{P}_{U,x}\subset\mc{S}_{U/k}$ the packet of cuspidal sections associated with $x$. We have that $\mc{S}_{k(X)/k}\simeq\projlim_{U}\mc{S}_{U/k}$ by \cite[Lemma 259]{sti13}, and the cuspidal birational sections over $x$ identify naturally with $\projlim_{U}\mc{P}_{U,x}\subset\mc{S}_{k(X)/k}$. Since $\iota_{U}(s)\in\mc{P}_{U,x}$ for every $U\subset X$ and $s$ is the limit of the sections $\iota_{U}(s)$, we get that $s$ is cuspidal associated with $x$ too.
\end{proof}

\begin{TC}
	The following are equivalent.
	\begin{itemize}
		\item The birational section conjecture holds.
		\item For every number field $k$, every section $z\in\mc{S}_{k(\P_{1})/k}$ and every open subset $U\s\P^{1}$, there exists an open subset $V\s U$ and a $\pi_{1}$ section $r$ of $U\times V\setminus\Delta \to V$ such that the specialization $r_{v}$ is equal to $\iota_{U\setminus\{v\}}(z)$ for every $v\in V(k)$.
	\end{itemize}
	\begin{proof}
		%~ \autoref{hilbgen} implies that $(iii)$ and $(iv)$ are equivalent. 
            Thanks to \autoref{bireq}, we may use the alternative formulation of the birational section conjecture.
    
		Assume that the birational section conjecture holds and let $z$, $U$ be as in the statement. Up to passing to an open subset, we can assume that $U$ is hyperbolic. By hypothesis, $z$ is associated with a unique rational point $p\in\P^{1}(k)$. If $p\in U$, choose $V=U\setminus\{p\}$, the morphism $(p,\id_{V}):V\to U\times V\setminus\Delta$ defines the desired $\pi_{1}$ section. If $p\not \in U$, choose $V=U$. Since $\iota_{U}(z)_{k(t)}$ is cuspidal then it lifts to a cuspidal section $z'$ of $U_{k(t)}\setminus\{\delta\}$. Let $v\in V(k)$ be a rational point, since $z'$ is cuspidal it is b.l. and its specializations at $v$ are b.l. too by \autoref{spebir}. By hypothesis, this implies that all the specializations of $z'$ at $v$ are cuspidal liftings of $\iota_{U}(z)$, and hence coincide with $\iota_{U\setminus\{v\}}(z)$ since there is only one cuspidal lifting of $\iota_{U}(z)$ by \cite[Theorem 17 (2)]{sti12}. In particular, the specializing loop at $v$ is constant by \autoref{uniqueloop1}. By repeating the argument after base changing to every finite extension of $k$, we get that all the specializing loops of $z'$ at closed points of $U$ are constant. As a consequence, $z'$ extends to the desired section $U\to\Pi_{(U\times U\setminus\Delta)/U}$ by \cite[Corollary A.3]{fed}.
  
            Assume that the second condition holds. Thanks to \cite[Theorem C]{saty21}, it is enough to consider number fields. Let $k$ be a number field, $X$ a smooth projective curve over $k$, $s\in\mc{S}_{X/k}$ a b.l. section with lifting $z\in\mc{S}_{k(X)/k}$ and $\nu$ a place of $k$, we want to prove that $s$ is geometric or cuspidal. With an argument analogous to that of \hyperlink{TA}{Theorem A}, we can reduce to the case in which $X=U$ is an open subset of $\P^{1}$ and to proving that $x_{\nu}\in\P^{1}$ is $k$-rational.
		
		Let $V\s U$, $r\in\mc{S}_{(U\times V\setminus\Delta)/V}$ be as given by hypothesis. If $U'\s U$ is another open subset and $V',r'$ are given analogously, we can shrink $V'$ and assume that $V'\s V$, then \autoref{hilbsp} implies that $r'=r|_{V'}$. Passing to the limit along open subsets of $\P^{1}$ we get a Galois section of $\spec k(\A^{1})\otimes_{k}k(t)\setminus\{\delta\}$ which lifts $s_{k(t)}$, hence $x_{\nu}$ is $k$-rational by \autoref{simple}.	
	\end{proof}
\end{TC}

%~ \printbibliography{}
\bibliographystyle{amsalpha}
\bibliography{tbir}

\providecommand{\bysame}{\leavevmode\hbox to3em{\hrulefill}\thinspace}
\providecommand{\MR}{\relax\ifhmode\unskip\space\fi MR }
\providecommand{\MRhref}[2]{%
  \href{http://www.ams.org/mathscinet-getitem?mr=#1}{#2}
}
\providecommand{\href}[2]{#2}
\begin{thebibliography}{{Gro}97}

\bibitem[AGV08]{agv08}
D.~{Abramovich}, T.~{Graber}, and A.~{Vistoli}, \emph{{G}romov-{W}itten theory
  of {D}eligne-{M}umford stacks}, American Journal of Mathematics \textbf{130}
  (2008), no.~5, 1337--1398.

\bibitem[And74]{anderson}
Michael~P. Anderson, \emph{Exactness properties of profinite completion
  functors}, Topology \textbf{13} (1974), no.~3, 229--239.

\bibitem[Bre21a]{fed}
Giulio Bresciani, \emph{Essential dimension and pro-finite group schemes}, Ann.
  Sc. Norm. Super. Pisa Cl. Sci. (5) \textbf{22} (2021), no.~4, 1899--1936.
  \MR{4360607}

\bibitem[Bre21b]{anab}
\bysame, \emph{Some implications between {G}rothendieck's anabelian
  conjectures}, Algebr. Geom. \textbf{8} (2021), no.~2, 231--267. \MR{4174290}

\bibitem[Bre23]{scfg}
\bysame, \emph{On the section conjecture over fields of finite type},
  arxiv:1911.03234, 2023.

\bibitem[BV15]{bv15}
N.~{Borne} and A.~{Vistoli}, \emph{{The Nori fundamental gerbe of a fibered
  category}}, Journal of Algebraic Geometry \textbf{24} (2015), 311--353.

\bibitem[BV19]{bv19}
Niels Borne and Angelo Vistoli, \emph{Fundamental gerbes}, Algebra Number
  Theory \textbf{13} (2019), no.~3, 531--576.

\bibitem[BV23]{giulio-angelo-valuative}
Giulio Bresciani and Angelo Vistoli, \emph{An arithmetic valuative criterion
  for proper maps of tame algebraic stacks}, Manuscripta Math. (2023),
  \url{https://doi.org/10.1007/s00229--023--01491--6}.

\bibitem[EH08]{eh08}
H\'{e}l\`ene Esnault and Ph\`ung~H\^{o} Hai, \emph{Packets in {G}rothendieck's
  section conjecture}, Adv. Math. \textbf{218} (2008), no.~2, 395--416.

\bibitem[GJZZ08]{gjz08}
F.~Grunewald, A.~Jaikin-Zapirain, and P.~A. Zalesskii, \emph{Cohomological
  goodness and the profinite completion of {B}ianchi groups}, Duke Math. J.
  \textbf{144} (2008), no.~1, 53--72.

\bibitem[Gro71]{sga1}
Alexandre Grothendieck, \emph{Rev\^{e}tements \'{e}tales et groupe
  fondamental}, Lecture Notes in Mathematics, Vol. 224, Springer-Verlag,
  Berlin-New York, 1971, S\'{e}minaire de G\'{e}om\'{e}trie Alg\'{e}brique du
  Bois Marie 1960--1961 (SGA 1), Dirig\'{e} par Alexandre Grothendieck.
  Augment\'{e} de deux expos\'{e}s de M. Raynaud.

\bibitem[{Gro}97]{gro97}
A.~{Grothendieck}, \emph{Brief an {G}. {F}altings}, Geometric Galois actions,
  1, London Math. Soc. Lecture Note Ser., vol. 242, Cambridge University Press,
  1997, With an English translation on pp. 285--293, pp.~49--58.

\bibitem[Koe05]{koe05}
Jochen Koenigsmann, \emph{On the `section conjecture' in anabelian geometry},
  J. Reine Angew. Math. \textbf{588} (2005), 221--235.

\bibitem[Mat55]{mat55}
Arthur Mattuck, \emph{Abelian varieties over {$p$}-adic ground fields}, Ann. of
  Math. (2) \textbf{62} (1955), 92--119.

\bibitem[Ols16]{olsson}
Martin Olsson, \emph{Algebraic spaces and stacks}, American Mathematical
  Society Colloquium Publications, vol.~62, American Mathematical Society,
  Providence, RI, 2016.

\bibitem[Sch15]{sch15}
Johannes Schmidt, \emph{Homotopy rational points of {B}rauer-{S}everi
  varieties}, arxiv:1503.08108, 2015.

\bibitem[ST21]{saty21}
Mohamed Sa\"{i}di and Michael Tyler, \emph{On the birational section conjecture
  over finitely generated fields}, Algebra Number Theory \textbf{15} (2021),
  no.~2, 435--460. \MR{4243653}

\bibitem[Sti12]{sti12}
Jakob Stix, \emph{On cuspidal sections of algebraic fundamental groups},
  Galois-Teichm\"uller theory and arithmetic geometry. Selected papers based on
  the presentations at the workshop and conference, Kyoto, Japan, October
  25--30, 2010, Tokyo: Mathematical Society of Japan, 2012, pp.~519--563.

\bibitem[{Sti}13]{sti13}
J.~{Stix}, \emph{Rational points and arithmetic of fundamental groups}, Lecture
  Notes in Mathematics, no. 2054, Springer, 2013.

\bibitem[Sti15]{sti15}
Jakob Stix, \emph{On the birational section conjecture with local conditions},
  Invent. Math. \textbf{199} (2015), no.~1, 239--265.

\bibitem[Tam97]{tam97}
Akio Tamagawa, \emph{The {G}rothendieck conjecture for affine curves},
  Compositio Math. \textbf{109} (1997), no.~2, 135--194.

\bibitem[TV18]{tv18}
Mattia Talpo and Angelo Vistoli, \emph{Infinite root stacks and quasi-coherent
  sheaves on logarithmic schemes}, Proc. Lond. Math. Soc. (3) \textbf{116}
  (2018), no.~5, 1187--1243.

\bibitem[Wic12]{wic12}
Kirsten Wickelgren, \emph{{$n$}-nilpotent obstructions to {$\pi_1$} sections of
  {$\Bbb P^1-\{0,1,\infty\}$} and {M}assey products}, Galois-{T}eichm\"{u}ller
  theory and arithmetic geometry, Adv. Stud. Pure Math., vol.~63, Math. Soc.
  Japan, Tokyo, 2012, pp.~579--600.

\end{thebibliography}

\end{document}